\documentclass[journal,twoside]{IEEEtran}
\usepackage{cite}
\usepackage{amsmath,amssymb,amsfonts}
\usepackage{algorithmic}
\usepackage{graphicx}
\usepackage{algorithm,algorithmic}
\usepackage{hyperref}
\usepackage{textcomp}
\usepackage{mathrsfs}
\usepackage{epstopdf} 
\usepackage{algorithm}
\usepackage{amsthm}
\usepackage{algorithmic}
\usepackage{color,xcolor}
\usepackage{amsmath,amsfonts,amssymb}
\usepackage{array}
\usepackage{multirow} 
\usepackage{arydshln}
\usepackage{nicematrix}

\newcommand{\N}{\ensuremath{\mathbb{N}}}

\newcommand{\R}{\ensuremath{\mathbb{R}}}

\def \R {\mathbb{R}}
\def \p {\partial}

\def \N {\mathcal{N}}
\def \bx {{\bf x}}
\def \blambda {{\bf\lambda}}
\def \bk {{\bf k}}
\def \be {{\bf e}}

\def \bq {{\bf q}}
\def \bw {{\bf\omega}}
\def \bgamma {{\bf\gamma}}
\def \bdelta {{\bf\delta}}   
\def \bX {{\bf X}}

\usepackage{graphicx}          

\newtheorem{THEOREM}{Theorem}[section]

\newtheorem{theorem}[THEOREM]{Theorem}
\newtheorem{proposition}[THEOREM]{Proposition}

\newtheorem{remark}[THEOREM]{Remark}
\newtheorem*{assumption*}{Assumption}

\theoremstyle{remark}
\newtheorem{example}[THEOREM]{Example}
\def\BibTeX{{\rm B\kern-.05em{\sc i\kern-.025em b}\kern-.08em
    T\kern-.1667em\lower.7ex\hbox{E}\kern-.125emX}}
\begin{document}
\title{A decomposition method in the multivariate feedback particle filter via tensor product Hermite polynomials}
\author{Ruoyu Wang$^{1}$, {\it IEEE student member} and Xue Luo$^{2,\sharp}$, {\it IEEE senior member}
\thanks{*This work is financially supported by National Natural Science Foundation of China (Grant No. 12271019) and the National Key R\&D Program of China (Grant No. 2022YFA1005103).} 
\thanks{$^{1}$R. Wang is with School of Mathematical Sciences, Beihang University, Beijing, P. R. China 102206. 
        {\tt\small WangRY@buaa.edu.cn}}%
\thanks{$^{2}$X. Luo is with School of Mathematical Sciences and Key Laboratory of Mathematics, Informatics and Behavioral Semantics (LMIB), Beihang University, Beijing, P. R. China 100191.
        {\tt\small xluo@buaa.edu.cn}}%
\thanks{$^{\sharp}$X. Luo is the corresponding author.}
}

\maketitle

\begin{abstract}
The feedback particle filter (FPF), a resampling-free algorithm proposed over a decade ago, modifies the particle filter (PF) by incorporating a feedback structure. Each particle in FPF is regulated via a feedback gain function (lacking a closed-form expression), which solves a Poisson’s equation with a probability-weighted Laplacian, as derived in \cite{YMM:13, YLMM:16}. While approximate solutions to this equation have been extensively studied in recent literature, no efficient multivariate algorithm exists. In this paper, we focus on the decomposition method for multivariate gain functions in FPF, which has been proven efficient for scalar FPF with polynomial observation functions. Its core is splitting the Poisson’s equation into two exactly solvable sub-equations. Key challenges in extending it to multivariate FPF include ensuring the invertibility of the coefficient matrix in one sub-equation and constructing a weighted-radial solution in the other. The proposed method’s computational complexity grows at most polynomially with the state dimension, a dramatic improvement over the exponential growth of most particle-based algorithms. Numerical experiments compare the decomposition method with traditional methods: the extended Kalman filter (EKF), PF, and FPF with constant-gain or kernel-based gain approximations. Results show it outperforms PF and FPF with other gain approximations in both accuracy and efficiency, achieving the shortest CPU time among methods with comparable performance.
\end{abstract}

\begin{IEEEkeywords}
Nonlinear filtering; multivariate feedback particle filter; Hermite spectral method; a weighted-radial solution.
\end{IEEEkeywords}

\section{Introduction}
\label{sec:introduction}
\IEEEPARstart{N}{onlinear} filtering (NLF) is a complex discipline dedicated to extracting valuable information from data corrupted by noisy sensors. It is the cornerstone of numerous applications, covering target tracking and navigation, air traffic management, meteorological monitoring, geographic surveying, geophysical measurement, remote sensing, autonomous navigation, and robotics fields, as documented by \cite{BLK:01}.

The NLF problem of diffusion processes is modeled by
\begin{equation}\label{eqn-NLF}
\left\{\begin{aligned}
dX_t=&g(X_t)dt+\sigma(X_t)dB_t\\
dZ_t=&h(X_t)dt+dW_t
\end{aligned}\right.,
\end{equation}
where $X_t\in\R^d$ is the state, $Z_t\in\R^m$ is the observation, $\{B_t\},\{W_t\}$ are two mutually independent standard Wiener processes taking values in $\R^d$ and $\R^m$, respectively. The drift function $g(\cdot):\ \R^d\rightarrow\R^d$ and the diffusion function $\sigma(\cdot):\ \R^d\rightarrow\R^{d\times d}$ are Lipshitz continuous, and the observation function $h=(h_1,h_2,\cdots,h_m)^T:\ \R^d\rightarrow\R^m$, where $h_j$ is the $j$th component of the column vector $h$, are bounded continuous.

Prior to the 1990s, the predominant methodologies for filtering applications centered on Kalman filter (KF) frameworks \cite{K:60,KB:61} and their nonlinear extensions, including the extended Kalman filter (EKF) \cite{SS:66,JA:70}, ensemble Kalman filter (EnKF) \cite{E:94,EG:03}, and related variants. However, when confronted with nonlinear systems, these methodologies exhibit significant limitations stemming from both signal model nonlinearities and measurement model constraints. Such nonlinear characteristics often induce non-Gaussian multi-modal conditional distributions, a scenario where conventional KF and EKF implementations demonstrate suboptimal performance \cite{BSN:04}.

The advent of particle filter (PF) and its variances represented a paradigm shift in NLF theory, gaining substantial traction in engineering applications due to their capability to handle non-Gaussian distributions \cite{CASV:04,GSS:93}. Nevertheless, practical implementations of PFs are constrained by the requirement for extensive sample populations to achieve sufficient approximation accuracy of the posterior probability density. Operational challenges including sample impoverishment, the curse of dimensionality, and variance escalation have been demonstrated to undermine the robustness \cite{AMGC:02}.

In their pioneering work \cite{YMM:13,YLMM:16}, Yang et al. proposed an innovative method named the feedback particle filter (FPF), which integrates controlled dynamic behaviors into each particle. This framework combines two key concepts-feedback control architecture and mean field game theory-enabling the derivation of approximate solutions to NLF problems. Compared with the PF, FPF achieves two critical improvements:
\begin{enumerate}
    \item[1)] It eliminates the need for resampling. Thanks to its inherent feedback mechanism, FPF avoids the ``particle degeneracy" issue-a major flaw in conventional methods where particles lose diversity and become ineffective over time.
    \item[2)] It adjusts particle states dynamically via closed-loop control. This ensures stable estimation accuracy even when dealing with nonlinear system dynamics, which often disrupt the performance of traditional filters.
\end{enumerate}

Addressing the calculation of the gain function $K$ is central in the efficiency of the FPF. It has been numerically verified by Surace et al. \cite{SKP:19} that once the gain function is accurately determined, the performance of the FPF can surpass that of almost all traditional NLF algorithms. In the quest for an accurate and efficient approximate gain function for the FPF, several methods have been developed. Initial work by Yang et al. \cite{YLMM:13} introduced a constant-gain approximation, which performs a straightforward averaging of particle controls. While the Galerkin method \cite{YMM:13, YLMM:16} enhanced accuracy, its requirement for pre-defined basis functions results in poor scalability and Gibbs phenomena in high-dimensional settings \cite{TM:16}. To address these shortcomings, Berntorp et al. \cite{BG:16} employed proper orthogonal decomposition for adaptive basis selection. Taghvaei et al. \cite{TM:16} then introduced a fundamentally different strategy with their kernel-based method, which requires no basis functions and has been analytically examined in \cite{TMM:17,TMM:20}. A comparative evaluation of these FPF approximations was conducted by Berntorp \cite{B:18}. The reproducing kernel Hilbert space (RKHS) approach introduced in \cite{AS:18} provides another basis-free algorithm. Later, Radhakrishnan et al. \cite{RM:19} combined it with the differential temporal-difference learning technique for gain function approximation in an on-line setting. However, as mentioned in the conclusion of \cite{RM:19}, the success of this method for high dimensional NLF problems are yet to be demonstrated. The survey \cite{AP:23} covers controlled interacting particle systems via FPF derived from optimal transport theory for NLF and optimal control, along with discussions on algorithmic derivations, comparative applications, and future research directions. In our recent work \cite{WML:25}, we developed a decomposition algorithm for the FPF in one-dimensional settings. This approach addresses the associated Poisson's equation by decomposing it into two sub-equations that admit exact solutions, particularly when the observation function is polynomial. 

In this paper, we further investigate the implementability of the decomposition method-Proposition \ref{prop-2.3} (Proposition III.1, \cite{WML:25})-for multivariate FPF, aiming to derive more efficient and accurate approximate gain functions. As outlined in \cite{WML:25}, the core idea of the decomposition method is to split the Poisson's equation into two sub-equations \eqref{eqn-1}-\eqref{eqn-71}; these can be explicitly and easily solved for scalar NLF problems, as shown in \eqref{eqn-2.1}-\eqref{eqn-2.16}, Theorem \ref{thm-d=1}. However, the key distinction between scalar and multivariate cases lies in the explicit solvability of \eqref{eqn-3.2}-\eqref{eqn-3.1} in Proposition \ref{prop-2.3}, which was briefly mentioned in Section III.C, \cite{WML:25}. The Galerkin spectral method, using tensor product Hermite polynomials, provides a potential approach to solve \eqref{eqn-3.1} and determine the corresponding constants $C_j^i$, $i=1,\cdots,d$, $j=1,\cdots,m$. Notably, unlike the scalar case, the invertibility of the coefficient matrices preceding the tensor product Hermite polynomials requires proof. For \eqref{eqn-3.2}, constructing an exact solution satisfying the boundary constraint is more challenging than the direct integration used in the scalar case, due to the presence of the divergence operator $\nabla^T$. The major contribution of this paper is to resolve these two issues and develop an implementable multivariate FPF via the decomposition method, when the observation is polynomial of the states. The main result is summarized in Theorem \ref{thm-3.5}. 

The organization of this paper is as follows: Section \ref{sec-2} reviews the multivariate FPF and the decomposition method for the scalar case. Section \ref{sec-3} focuses on deriving the exact solutions to the two sub-equations of the decomposition method, which are pertinent to the multivariate FPF. Section \ref{sec-4} presents the results of numerical experiments, including those on the ship tracking problem and the Lorenz oscillator system. The conclusion is drawn in the end.

\section{Preliminaries}\label{sec-2}

\subsection{Multivariate feedback particle filter (FPF)}\label{sec-2.1}

\IEEEPARstart{I}{n} FPF \cite{YLMM:16}, the $i$-th particle is governed by a controlled system
\begin{equation}\label{eqn-2.20}
dX_t^i=g(X_t^i)dt+\sigma(X_t^i)dB_t^i+dU_t^i,
\end{equation} 
$i=1,\cdots,N_p$, where 
\begin{align}\label{eqn-2.3}
    dU_t^i=u(X_t^i,t)dt+K(X_t^i,t)dZ_t
\end{align}
    is the $i$-th particle's control, $X_t^i$ is the $i$-th particle's state at time $t$, and $\{B_t^i\}$ are mutually independent standard Wiener processes. The optimal control $K$ is the minimizer of an optimization problem. Its Euler-Lagrange boundary value problem is obtained via the analysis of first variation. The gain function $K:\R^d\times\R_+\rightarrow\R^{d\times m}$ is the solution to 
\begin{equation}\label{eqn-BVP}
	\nabla^T(p_tK)=-\left(h-p_t[h]\right)^Tp_t,
\end{equation}
and $u:\R^d\times \R_+\rightarrow\R^d$ is obtained by
\begin{equation}\label{eqn-u}
u=-\frac12K\left(h+p_t[h]\right)+\Omega(\bx,t),
\end{equation}
with $p_t[h]:=\int_{\R^d}h(\bx)p_t(\bx)d\bx$, where $\Omega=(\Omega_1,\cdots,\Omega_d)^T$ is the Wong-Zakai correction term
\begin{equation}\label{eqn-Omega}
    \Omega_l(\bx,t) :=\frac12\sum_{k=1}^d\sum_{s=1}^mK_{ks}(\bx,t)\frac{\partial K_{ls}}{\partial x_k}(\bx,t),
\end{equation}
$l=1,\cdots,m$. Yang et al. showed that with the optimal control pair $(K,u)$ in \eqref{eqn-BVP}-\eqref{eqn-u} the conditional density of the particles, denoted as $p_t$, matches that of the true state $X_t$, under some mild conditions, see Theorem 1, \cite{YLMM:16}.

The evolution of the density $p_t$ satisfies the forward Kolmogorov equation, Proposition 3.2.1, \cite{Y:14}:
\begin{align}\label{eqn-p_t}\notag
    dp_t  =& \mathcal{L}^*p_tdt-\nabla^T(p_tK)dZ_t-\nabla^T(p_tu)dt\\
    &+\frac12\sum_{l,k=1}^d\frac{\p^2}{\p X_l\p X_k}(p_t[KK^T]_{lk})dt.
\end{align}
It is suggested by \cite{Y:14} that $p_t$ can be approximated by the empirical distribution of the controlled particles $\{X_t^i\}_{i=1}^{N_p}$, i.e.
\begin{equation}
    p_t^{(N_p)}(A)=\frac1{N_p}\sum_{i=1}^{N_p}1_{\{F(X_0^i,B_{[0,t]}^i;Z_{[0,t]})\in A\}}, 
\end{equation}
where $X_t^i=F(X_0^i,B_{[0,t]}^i;Z_{[0,t]})$ is a functional representation for each $i$, with the notation $Z_{[0,t]}$ signifies the entire observation path $\{Z_s:0\leq s\leq t\}$, $F$ is a continuous functional of the sample path $\{B_{[0,t]}^i,Z_{[0,t]}\}$ along with the initial condition $X_0^i$, \cite{K:90}. The almost sure convergence of $p_t^{(N_p)}$ to $p_t$ is guaranteed by the Law of Large Numbers. Moreover, we approximate the empirical distribution $p_t^{(N_p)}$ by the Gaussian mixture 
  \begin{equation}\label{eqn-2.4}
    p_t^{N_p,\Sigma}(\bx)=\frac1{N_p}\sum_{i=1}^{N_p}\N(\bx;X^i_t,\Sigma),
    \end{equation}
    where
    \begin{align}\label{eqn-Gaussian multi}
	\mathcal{N}(\bx;X_t^i,\Sigma)
    =&\frac{1}{(2\pi)^{\frac{d}{2}}|\Sigma|^{\frac{1}{2}}}\exp\left\{-\frac{1}{2}(\bx-X_t^i)^T\Sigma^{-1}(\bx-X_t^i)\right\}
  \end{align}
  is the multivariate Gaussian density function. Thus, \eqref{eqn-BVP} is approximated by
    \begin{equation}\label{eqn-pN-K}
    \nabla^T\left(p_t^{N_p,\Sigma}K\right)=-\left(h-p_t^{N_p,\Sigma}[h]\right)^Tp_t^{N_p,\Sigma},
    \end{equation}
where
   \begin{equation}\label{eqn-hat hN}
    p_t^{N_p,\Sigma}[h]=\int_{\R^d}h(\bx)p^{N_p,\Sigma}_t(\bx)d\bx.
    \end{equation}

In the sequel, we shall develop a decomposition method to solve \eqref{eqn-pN-K} explicitly, when the observation function $h(\bx)$ is a polynomial in $\bx$.

\subsection{Decomposition method in the scalar case}\label{sec-2.2}

From Section \ref{sec-2.1}, it is not hard to see that finding the exact/approximate solution to \eqref{eqn-pN-K} is crucial to guarantee the accuracy and efficiency of the FPF. 

In the scalar case, i.e. $d=1$, 
\eqref{eqn-pN-K} can be solved by direct integration
\begin{align}\label{eqn-direct K}
	K(x,t)=\frac1{p_t^{N_p,\varepsilon}(x)}\int_{-\infty}^x-\left[h(y)-p_t^{N_p,\varepsilon}[h]\right]p_t^{N_p,\varepsilon}(y)dy,
\end{align}
satisfying the boundary condition
\begin{equation}\label{eqn-90}
	\lim_{x\rightarrow-\infty}K(x,t)p^{N_p,\varepsilon}_t(x)=0,
\end{equation}
where the superscript $\Sigma$ in \eqref{eqn-pN-K} degenerates to $\varepsilon$. 

On the contrary, in the multivariate case, the direct integration can't work at all, due to the divergence operator $\nabla^T$ in \eqref{eqn-pN-K}. Thus, it is necessary to find a feasible method. In \cite{WML:25}, the authors of this paper proposed a decomposition method to obtain the gain function when $d=1$ under the boundary condition \eqref{eqn-90}, when the observation $h(x)$ is a polynomial in $x\in\R$.

Let us briefly recall the strategy. First, it is observed that \eqref{eqn-pN-K} can be decomposed into two components.

\begin{proposition}[Corollary III.2, \cite{WML:25}]\label{prop-2.1}
The gain function $K(x)$ in \eqref{eqn-pN-K} is given by
\begin{align}\label{eqn-50}
	K(x)
    =&\frac1{\sum_{i=1}^{N_p}\mathcal{N}(x;X_t^i,\varepsilon)}\sum_{i=1}^{N_p}\left[\mathcal{N}(x;X_t^i,\varepsilon)K^i(x)+K^i_0(x)\right],
\end{align}
where the functions $K^i(x)$ and $K_0^i(x)$  satisfy
\begin{equation}\label{eqn-1}
	\left[\N(x;X_t^i,\varepsilon){K^i}(x)\right]'=-\left(h(x)-C^i\right)\N(x;X_t^i,\varepsilon),
\end{equation}
 and 
\begin{equation}\label{eqn-71}
	{K_0^i}'(x)=\left(p_t^{N_p,\varepsilon}[h]-C^i\right)\N(x;X_t^i,\varepsilon),
\end{equation}
for any constant $C^i$, respectively.
\end{proposition}

Then, under the assumption that $h(x)$ is a polynomial in $x$, one proceeds 
\begin{enumerate}
    \item[ 1)] Equation \eqref{eqn-1} is solved exactly by the Galerkin Hermite spectral method and the backward recursion to determine $C^i$. 
    \item[2)] With the $C^i$ in 1), equation \eqref{eqn-71} is integrated directly to obtain the exact solution.
\end{enumerate}
By the procedure above, the exact gain function for $d=1$ has been obtained:
\begin{theorem}[Theorem III.3, \cite{WML:25}]\label{thm-d=1}
	When $d=1$ and $\displaystyle h(x)=\sum_{k=0}^pa_kH_{k}(x)$ is a polynomial of degree $p$, for some given $p\geq1$, where $H_k(x)$ represents the Hermite polynomial of degree $k$. The exact solution to \eqref{eqn-pN-K} is given by \eqref{eqn-50}, where
\begin{align}\label{eqn-2.1}
	K^i(x)=&\sum_{l=0}^{p-1}\tilde K^i_lH_l(x),
\end{align}
 where the coefficients $\{\tilde K^i_l\}_{l=0}^{p-1}$ are calculated through backward recursion
\begin{equation}\label{eqn-2.15}
	\tilde K_k^i=2\varepsilon a_{k+1}+2(\varepsilon-1)(k+2)\tilde K_{k+2}^i+2X_t^i\tilde K_{k+1}^i,
\end{equation}
for $k=p-1,\cdots,1,0$, with the initial values $\tilde K^i_{p+1}=\tilde K^i_{p}\equiv0$, and
\begin{equation}\label{eqn-2.16}
	C^i:=a_0+\frac{X_t^i}{\varepsilon}\tilde K^i_0+\left(2-\frac1{\varepsilon}\right)\tilde K^i_1,
\end{equation}
and
\begin{align}\label{eqn-2.2}
       K_0^i(x)=\frac{p_t^{N_p,\varepsilon}[h]-C^i}2\textup{erf}\left(\frac{x-X_t^i}{\sqrt{2\varepsilon}}\right),
\end{align}
with $\textup{erf}(x)=\frac2{\sqrt\pi}\int_0^xe^{-\eta^2}d\eta$ being the error function.
\end{theorem}

The brief discussions on the applicability to the multivariate FPF are mentioned in Section III.C, \cite{WML:25}. It has already shown that the multivariate version of Proposition \ref{prop-2.1} holds. But the exact solvability of \eqref{eqn-3.2}-\eqref{eqn-3.1} are left to be addressed.

\begin{proposition}[Proposition III.1, \cite{WML:25}]\label{prop-2.3}
    For each $j=1,\cdots,m$, the $j$-th column of the gain function $K_{\cdot j}=\nabla\varphi_j(\bx)$ is given by
\begin{align}\label{eqn-3.10}\notag
	\nabla\varphi_j(\bx)
	=&\frac{1}{\sum_{i=1}^{N_p}\mathcal{N}(\bx;X_t^i,\Sigma_i)}\\
    &\cdot\sum_{i=1}^{N_p}\left[\nabla_{\blambda^i}\psi_j^i(\bx)+\mathcal{N}(\bx;X_t^i,\Sigma_i)\nabla\varphi_j^i(\bx)\right],
\end{align}
with $\blambda^i=(\lambda_1^i,\cdots,\lambda_d^i)$ being the eigenvalues of $\Sigma_i^{-1}$, and $\nabla_{{\blambda}^i}\psi_j^i(\bx):=\left(\frac{1}{\lambda_l^i}\frac{\partial\psi_j^i(\bx)}{\partial x_l}\right)_{l=1}^{d}$. The functions $\psi_j^i(\bx)$ and $\varphi_j^i(\bx)$ satisfy
\begin{equation}\label{eqn-3.2}
	\nabla^T\left(\nabla_{\mathbf{\lambda}^i}\psi_j^i(\bx)\right)=\left(p_t^{N_p,\Sigma}[h_j]-C_j^i\right)\mathcal{N}(\bx;X_t^i,\Sigma_i),
\end{equation}
and
	\begin{equation}\label{eqn-3.1}
		\nabla^T\left(\mathcal{N}(\bx;X_t^i,\Sigma_i)\nabla\varphi_j^i(\bx)\right)=\left(-h_j(\bx)+C_j^i\right)\mathcal{N}(\bx;X_t^i,\Sigma_i),
	\end{equation}
for any constant $C_j^i$, respectively.
\end{proposition}

Following the same strategy for $d=1$, the main difficulties in the multivariate case lie in both two steps:
\begin{enumerate}
    \item the implementability of the Galerkin spectral method and backward recursion in solving \eqref{eqn-3.1};
    \item the construction of the exact solution to \eqref{eqn-3.2}, where the direct integration is not applicable anymore.
\end{enumerate}

In Section \ref{sec-3}, we shall focus on these two issues to give an explicit expression of the gain function satisfying certain conditions, like \eqref{eqn-50} in Proposition \ref{prop-2.1} and \eqref{eqn-2.1}-\eqref{eqn-2.2} in Theorem \ref{thm-d=1} for $d=1$.

\section{Decomposition Method for the multivariate FPF}\label{sec-3}

\IEEEPARstart{A}{s} suggested in \cite{YLMM:16}, let $\varphi(\bx)=(\varphi_{1}(\bx),\cdots,\varphi_m(\bx))$ be a $m$-vector valued function, such that the $j$-th column of the gain function $K_{\cdot j}(\bx)=\nabla\varphi_j(\bx)$, $j=1,\cdots,m$. Then \eqref{eqn-pN-K} becomes
\begin{align}\label{eqn-pN_phi}
	\nabla^T\left[p_t^{N_p,\Sigma}(\bx)\nabla\varphi_j(\bx)\right]=-\left(h_j(\bx)-p_t^{N_p,\Sigma}[h_j]\right)p_t^{N_p,\Sigma}(\bx).
\end{align}
Due to Proposition \ref{prop-2.3}, we shall obtain the {\it exact} solutions to \eqref{eqn-3.2}-\eqref{eqn-3.1}, when $h_j(\bx)$ are polynomials of at most degree $p$:
\begin{assumption*}
   \begin{equation}\label{eqn-3.14}
	h_j(\bx)=\sum_{\bk\in\Omega}a_{j,\bk}H_{\bk}(\bx),
	\end{equation}
for all $j=1,\cdots,m$, where $\Omega:=\{\bk\in\mathbb{N}^d:\,0\leq|\bk|_1\leq p\}$, $|\bk|_1=\sum_{l=1}^dk_l$, and $H_{\bk}(\bx)=\prod_{l=1}^d H_{k_l}(x_l)$ are the tensor product Hermite polynomials.
\end{assumption*}

\subsection{Decomposition Method }\label{DM}
As we pointed out before, the direct integration is not applicable to solve \eqref{eqn-pN_phi} when $d\geq2$, due to the divergence operator $\nabla^T$. In the following two subsections, we shall give an analytic solution to  \eqref{eqn-3.2} and propose an algorithm to yield an exact solution of \eqref{eqn-3.1} in the linear subspace spanned by the tensor product Hermite polynomials.

\subsubsection{A radially symmetric solution to \eqref{eqn-3.2}}

We shall give an explicit expression of a radially symmetric solution to \eqref{eqn-3.2}.
\begin{proposition}[Radially symmetric solution]\label{P1} 
Equation \eqref{eqn-3.2} has a radially symmetric solution
\begin{align}\label{eqn-15}
	\nabla_{\blambda^i}\psi_j^i(\bx)
=(\bx-X_t^i)\frac{\left(p_t^{N_p,\Sigma}[h_j]-C_j^i\right)}{2\pi^{\frac{d}{2}}|\Sigma_i|^{\frac12}}\gamma\left(\frac d2,\frac{r^2}{2}\right)r^{-d},
\end{align}
with the weighted radial $r:=\sqrt{\sum_{l=1}^{d}\lambda_l^i(x_l-X_{t,l}^i)^2}$, $\left\{\lambda_l^i\right\}_{l=1}^d$ being the eigenvalues of $\Sigma_i^{-1}$, such that $\nabla_{\lambda^i}\psi_j^i(X_t^i)={\bf 0}$, for $j=1,\cdots,m$, $i=1,\cdots,N_p$, where $\gamma(s,x):=\int_0^xt^{s-1}e^{-t}dt$ is the lower incomplete Gamma function.
\end{proposition}
\begin{remark}\label{rmk-10}
	For $d=1$, the solution to \eqref{eqn-3.2} is 
\begin{align}\label{eqn-58}\notag
	\nabla_{\lambda^i}\psi^i(x)=&(x-X_t^i)\frac{\left(p_t^{N_p,\Sigma}[h]-C^i\right)\sqrt{\lambda^i}}{2\sqrt\pi r}\gamma\left(\frac12,\frac{r^2}2\right)\\\notag
    =&\frac{x-X_t^i}{|x-X_t^i|}\frac{p_t^{N_p,\Sigma}[h]-C^i}2\textup{erf}\left(\frac{|x-X_t^i|}{\sqrt{2\varepsilon_i}}\right)\\
	=&\frac{p_t^{N_p,\Sigma}[h]-C^i}2\textup{erf}\left(\frac{x-X_t^i}{\sqrt{2\varepsilon_i}}\right),
\end{align}
since $\gamma\left(\frac12,x\right)=\sqrt\pi\textup{erf}\left(\sqrt x\right)$, $x>0$. This is exactly the solution $K_0^i(x)$ of \eqref{eqn-2.2} in Theorem \ref{thm-d=1} (Theorem III.3, \cite{WML:25}). It confirms that the decomposition method is an alternative way to solve \eqref{eqn-3.2} when $d=1$, as well as a feasible way to be extended to the multivariate case.
\end{remark}
\begin{IEEEproof}[Proof of Proposition \ref{P1}.]
  We first claim that the right-hand side of \eqref{eqn-3.2} is a function of $r$. Indeed, by the unit orthogonal decomposition of $\Sigma_i^{-1}=P_i^TD_iP_i$, with $P_i$ being the unit orthogonal matrix, $D_i=\textup{diag}\left(\lambda_1^i,\cdots,\lambda_d^i\right)$, the right-hand side of \eqref{eqn-3.2}, up to a constant $\left(p_t^{N_p,\Sigma}[h_j]-C_j^i\right)$, can be written as
\begin{align}\label{eqn-3.101}\notag
	&\frac{1}{(2\pi)^{\frac{d}{2}}|\Sigma_i|^{\frac{1}{2}}}\exp\left[-\frac{1}{2}(\bx-X_t^i)^T\Sigma_i^{-1}(\bx-X_t^i)\right]\\\notag
	=&\frac{\Pi_{l=1}^d\sqrt{\lambda_l^i}}{(2\pi)^{\frac{d}{2}}}\exp\left[-\frac{1}{2}(\bx-X_t^i)^TP_i^TD_iP_i(\bx-X_t^i)\right]\\\notag
	=&\frac{\Pi_{l=1}^d\sqrt{\lambda_l^i}}{(2\pi)^{\frac{d}{2}}}\exp\left[-\frac{1}{2}\sum_{l=1}^d\lambda_l^i(x_l-X_{t,l}^i)^2\right]\\
	=&\frac{\Pi_{l=1}^d\sqrt{\lambda_l^i}}{(2\pi)^{\frac{d}{2}}}\exp\left(-\frac12r^2\right),
\end{align}
where the second equality follows from the fact that \begin{align*}
    &(\bx-X_t^i)^TP_i^TD_iP_i(\bx-X_t^i)=\left|\sqrt{D_i}P_i(\bx-X_t^i)\right|^2\\
    =&\left|P_i\sqrt{D_i}(\bx-X_t^i)\right|^2=\left|\sqrt{D_i}(\bx-X_t^i)\right|^2\\
    =&\sum_{l=1}^d\lambda_l^i(x_l-X_{t,l}^i)^2.
\end{align*}

Based on the above observation, we shall look for a weighted radially symmetric solution to \eqref{eqn-3.2}, denoted as $\psi_j^i(r)$. Thus, the partial differential equation (PDE) \eqref{eqn-3.2} is reduced to an ordinary differential equation (ODE) with respect to $r$, that is,
  \begin{align}\label{eqn-3.5}\notag
	&\left(p_t^{N_p,\Sigma}[h_j]-C_j^i\right)\frac{\Pi_{l=1}^d\sqrt{\lambda_l^i}}{(2\pi)^{\frac{d}{2}}}\exp\left(-\frac12r^2\right)\\\notag
	\overset{\eqref{eqn-3.2},\eqref{eqn-3.101}}=&\nabla^T\left(\nabla_{\blambda^i}\psi_j^i(\bx)\right)\\\notag
 	=&\sum_{l=1}^d{\psi_j^i(r)}''\frac{\lambda_l^i(x_l-X_{t,l}^i)^2}{r^2}\\\notag
    &+\sum_{l=1}^d{\psi_j^i(r)}'\left(\frac{1}{r}-\frac{\lambda_l^i(x_l-X_{t,l}^i)^2}{r^3}\right)\\
	 =&{\psi_j^i}''(r)+\frac{d-1}{r}{\psi_j^i}'(r).
\end{align}
\begin{enumerate}
	\item When $d=1$, the first-order derivative term on the right-hand side of \eqref{eqn-3.5} vanishes and \eqref{eqn-3.5} is reduced to 
\begin{align}\label{eqn-3.102}
	{\psi^i}''(r)=\left(p_t^{N_p,\Sigma}[h]-C^i\right)\frac{\sqrt{\lambda^i}}{\sqrt{2\pi}}\exp\left(-\frac12r^2\right),
\end{align} 
with $r=\sqrt{\lambda^i}|x-X_t^i|$ and $\Sigma_i=\frac1{\lambda^i}$. Equation \eqref{eqn-3.102} can be solved by direct integration with respect to $r$, i.e. 
\begin{align}\label{eqn-57}\notag
	&{\psi^i}'(r)\\\notag
     =&{\psi^i}'(0)+\left(p_t^{N_p,\Sigma}[h]-C^i\right)\sqrt{\lambda^i}\left[\frac1{\sqrt{2\pi}}\int_{0}^re^{-\frac{s^2}2}ds\right]\\
	=&\frac{\left(p_t^{N_p,\Sigma}[h]-C^i\right)\sqrt{\lambda^i}}{2}\textup{erf}\left(\frac r{\sqrt{2}}\right),
\end{align}
with ${\psi^i}'(0)=0$. 
	\item When $d\geq2$, equation \eqref{eqn-3.5} is a first-order ODE of ${\psi_j^i}'(r)$. It is not hard to verify that $r^{1-d}$ is the general solution to the homogeneous part, and the special solution of the inhomogeneous part can be obtained by the method of variation of parameters. Indeed, letting ${\psi_j^{i*}}'(r)=\alpha(r)r^{1-d}$, with $\alpha(r)$ to be determined. Substituting ${\psi_j^{i*}}'(r)$ back to \eqref{eqn-3.5}, one has
\begin{align*}
	&{\psi_j^{i*}}''(r)+\frac{d-1}r{\psi_j^{i*}}'(r)
    =\alpha'(r)r^{1-d}\\
	=&\left(p_t^{N_p,\Sigma}[h_j]-C_j^i\right)\frac{\Pi_{l=1}^d\sqrt{\lambda_l^i}}{(2\pi)^{\frac{d}{2}}}\exp\left(-\frac12r^2\right).
\end{align*}
Thus, $\alpha(r)$ is solved by integrating from $0$ to $r$, i.e.
\begin{align*}
	\alpha(r)=&\alpha(0)+\frac{\left(p_t^{N_p,\Sigma}[h_j]-C_j^i\right)\Pi_{l=1}^d\sqrt{\lambda_l^i}}{(2\pi)^{\frac{d}{2}}}\\
    &\phantom{\left(\hat{h}_j-\right)}\cdot\int_0^rs^{d-1}\exp\left(-\frac12s^2\right)ds\\
	=&\alpha(0)+\frac{\left(p_t^{N_p,\Sigma}[h_j]-C_j^i\right)\Pi_{l=1}^d\sqrt{\lambda_l^i}}{2\pi^{\frac{d}{2}}}\gamma\left(\frac d2,\frac{r^2}2\right),
\end{align*}
where $\gamma(s,x)$ is the lower incomplete Gamma function. Therefore, the solution to \eqref{eqn-3.5} is 
\begin{align}\label{eqn-3.104}
	{\psi_j^i}'(r)=Cr^{1-d}+\frac{p_t^{N_p,\Sigma}[h_j]-C_j^i}{2\pi^{\frac{d}{2}}|\Sigma_i|^{\frac{1}{2}}}\gamma\left(\frac d2,\frac{r^2}2\right)r^{1-d},
\end{align}
where $C$ is an arbitrary constant. To uniquely determine this constant $C$, let us check the behavior of ${\psi_j^i}'(r)$ as $r\rightarrow0$. With the property $\frac{\gamma(s,r)}{r^s}\rightarrow\frac1s$, as $r\rightarrow0$, we have
$\frac{\gamma(\frac d2,\frac {r^2}2)}{r^d}\rightarrow\frac{2^{1-\frac d2}}d$. That is, the second term on the right-hand side of \eqref{eqn-3.104} is of order $r$, as $r\rightarrow0$, while the first term is of order $r^{1-d}\rightarrow\infty$. The constant $C$ is forced to vanish, i.e. ${\psi_j^i}'(0)=0$. That is,
\begin{equation}\label{eqn-56}
	{\psi_j^i}'(r)=\frac{p_t^{N_p,\Sigma}[h_j]-C_j^i}{2\pi^{\frac{d}{2}}|\Sigma_i|^{\frac{1}{2}}}\gamma\left(\frac d2,\frac{r^2}2\right)r^{1-d}.
\end{equation}
 \end{enumerate} 
By Remark \ref{rmk-10}, equation \eqref{eqn-57} is the same as \eqref{eqn-56} when $d=1$. Consequently, equation \eqref{eqn-15} follows immediately from 
$\displaystyle\nabla_{\blambda^i}\psi_j^i(\bx)
:=\left(\frac1{\lambda_l^i}\frac{\p\psi_j^i(\bx)}{\p X_l}\right)_{l=1}^d=\frac{\bx-X_t^i}r{\psi_j^i}'(r)$ and $\displaystyle\nabla_{\lambda^i}\psi_j^i(X_t^i)=\lim_{r\rightarrow0}(\bx-X_t^i)\frac{{\psi_j^i}'(r)}r={\bf 0}$.
\end{IEEEproof}

\subsubsection{An exact solution to \eqref{eqn-3.1}}

Recall the backward recursion of the coefficients in the scalar case illustrated in Fig. 1, Section III.B, \cite{WML:25}. We shall follow the same procedure for the multivariate case. By a direct computation, one sees that \eqref{eqn-3.1} is equivalent to 
\begin{align}\label{eqn-3.3}
	-(\bx-X_t^i)^T\Sigma_i^{-1}\nabla\varphi^i_j(\bx)+\bigtriangleup\varphi_j^i(\bx)=-h_j(\bx)+C_j^i,
\end{align}
for $j=1,\cdots,m$. By the Galerkin spectral method, let
\begin{equation}\label{eqn-3.4}
	\varphi^i_j(\bx)=\sum_{\bk\in\Omega}\tilde{\varphi}^i_{j,\bk}H_{\bk}(\bx),
\end{equation}
 where $\tilde{\varphi}^i_{j,\bk}$ are the coefficients to be determined by \eqref{eqn-3.3}, $\Omega:=\left\{\bk=(k_1,\cdots,k_d)\in\mathbb{N}^d:\,0\leq |\bk|_1\leq p\right\}$ is the index set and $H_{\bk}(\bx):=\Pi_{l=1}^d H_{k_l}(x_l)$ are the $d$-dimensional tensor product Hermite polynomials.

\begin{proposition}[Backward recursion]\label{P2}
   Under the assumption \eqref{eqn-3.14}, the coefficients $\tilde\varphi_{j,\bq}^i$, $1\leq|\bq|_1\leq p$, can be obtained by the backward recursion
\begin{align}\label{eqn-3.17}\notag
	&\sum_{l=1}^{d}\tilde{\varphi}_{j,\bq}^{i}\left(\Sigma_i^{-1}\right)_{ll}q_l+\sum_{\substack{l,m=1\\\notag l\not=m}}^d\tilde{\varphi}_{j,\bq+\be_m-\be_l}^i\left(\Sigma_i^{-1}\right)_{lm}(q_m+1)\\\notag
	=&a_{j,\bq}+\sum_{l,m=1}^d\tilde{\varphi}_{j,\bq+\be_m}^i\left(\Sigma_i^{-1}\right)_{lm}2(q_m+1)X_{t,l}^i\\
	&-\sum_{\substack{l,m=1\\\notag l\not=m}}^d\tilde{\varphi}_{j,\bq+\be_m+\be_l}^i\left(\Sigma_i^{-1}\right)_{lm}2(q_m+1)(q_l+1)\\
	&+\sum_{l=1}^d\tilde{\varphi}_{j,\bq+2\be_l}^i2\left[2-\left(\Sigma_i^{-1}\right)_{ll}\right](q_l+1)(q_l+2),
	\end{align}
with the initial condition $\tilde\varphi_{j,\bq}^i\equiv0$, for $|\bq|_1>p$ and the convention $\tilde\varphi_{j,\bq}^i\equiv0$, if some component in $\bq$ is negative, say $q_l<0$. Then $\varphi^i_j(\bx)$ in \eqref{eqn-3.4} is the solution of \eqref{eqn-3.3} with
    \begin{align}\label{eqn-Cji}
        C_j^i
     =&a_{j,{\bf0}}+\sum_{l,m=1}^d\tilde{\varphi}_{j,\be_m}^i\left(\Sigma_i^{-1}\right)_{lm}2X_{t,l}^i\\\notag
    &-\sum_{\substack{l,m=1,\\l\not=m}}^d\tilde{\varphi}_{j,\be_m+\be_l}^i2\left(\Sigma_i^{-1}\right)_{lm}
    +\sum_{l=1}^d\tilde{\varphi}_{j,2\be_l}^i4\left[2-\left(\Sigma_i^{-1}\right)_{ll}\right].
    \end{align} 
The notation $\be_l:=(0,\cdots,1,\cdots,0)$ is a zero vector except with only the $l$th component being $1$, and $X_{t,l}^i$ denotes the $l$th component of $X_t^i$.
 \end{proposition}

 \begin{figure}[th!]
    \centerline{\includegraphics[width=\columnwidth]{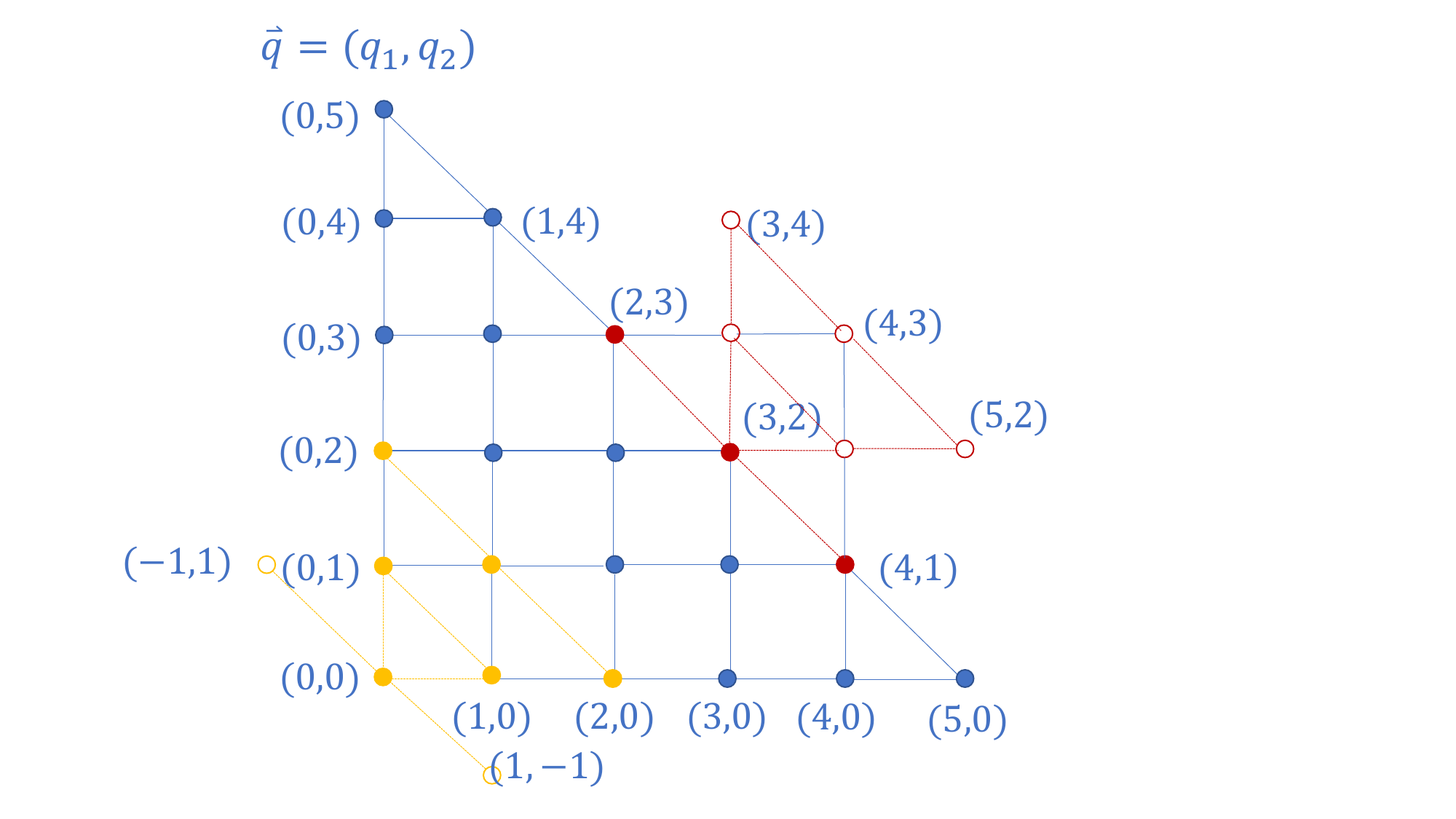}}
    \caption{The index set $\Omega$ is plotted to show the procedure of backward recursion in \eqref{eqn-3.17}-\eqref{eqn-Cji}, when $d=2$, $p=5$.}\label{fig-d2}
\end{figure}

To describe the backward recursion vividly, we present a concrete example.
\begin{example}
Let us take $d=2$, $p=5$. In Fig. \ref{fig-d2}, the index set $\Omega$ is plotted, using the circles at the coordinate $\bk=(k_1,k_2)$ to present $\tilde\varphi_{j,\bk}^i\equiv0$. From \eqref{eqn-3.17}, starting from $\tilde\varphi_{j,\bk}^i\equiv0$ with $|\bk|_1=6$ and $7$, the coefficients on the line $|\bk|_1=5$ are determined. Backward recursively, the coefficients $\tilde\varphi_{j,\bk}^i$ on each oblique line are determined till $|\bk|_1=1$. At last, the constant $C_j^i$ is determined by $\tilde\varphi_{j,\bk}^i$, $|\bk|_1=1,2$.
\end{example}

\begin{IEEEproof}[Proof of Proposition \ref{P2}]
By plugging \eqref{eqn-3.4} into \eqref{eqn-3.3}, one has 
\begin{align*}  
	&-h_j(\bx)+C^i_j	\\
	=&-(\bx-X_t^i)^T\Sigma_i^{-1}
\left[
\sum_{\bk\in\Omega}\tilde{\varphi}_{j,\bk}^{i}
\left(
\begin{array}{ccc} 2k_1H_{\bk-\be_1}(\bx)\\\notag
2k_2H_{\bk-\be_2}(\bx)\\
\vdots\\
2k_dH_{\bk-\be_d}(\bx)
\end{array}
\right)
\right]\\\notag
	&\phantom{2k_1H_{\bk-\be_1}(\bx)}+\sum_{\bk\in\Omega}\tilde{\varphi}_{j,\bk}^{i}\sum_{l=1}^{d}4k_l(k_l-1)H_{\bk-2\be_l}(\bx)\\\notag
    =&\sum_{\bk\in\Omega}\tilde{\varphi}_{j,\bk}^{i}\left[-\sum_{l,m=1}^d(x_l-X_{t,l}^i)\left(\Sigma_i^{-1}\right)_{lm}2k_mH_{\bk-\be_m}(\bx)\right.\\\notag
    &\phantom{\tilde{\varphi}_{j,\bk}^{i}\left[-\sum_{l,m=1}^d(x_l-x_{i,l})\right.}\left.+\sum_{l=1}^{d}4k_l(k_l-1)H_{\bk-2\be_l}(\bx)\right]\\\notag
    =&\sum_{\bk\in\Omega}\tilde{\varphi}_{j,\bk}^{i}\left\{-\sum_{\substack{l,m=1\\ l\not=m}}^{d}\left(\Sigma_i^{-1}\right)_{lm}2k_m\left[\frac{1}{2}H_{\bk-\be_m+\be_l}(\bx)\right.\right.\\\notag
    &\phantom{\sum_{\bk\in\Omega}\tilde{\varphi}_{j,\bk}^{i}\left\{-\sum_{\substack{l,m=1\\ l\not=m}}^{d}\left(\Sigma_i^{-1}\right)_{lm}2k_m\right.}\left.\left.+k_lH_{\bk-\be_m-\be_l}(\bx)\right]\right.\\\notag
	-&\left.\sum_{l=1}^{d}\left(\Sigma_i^{-1}\right)_{ll}2k_l\left[\frac{1}{2}H_{\bk}(\bx)+(k_l-1)H_{\bk-2\be_l}(\bx)\right]\right.\\\notag
	+&\left.\sum_{l,m=1}^{d}\left(\Sigma_i^{-1}\right)_{lm}2k_mX_{t,l}^iH_{\bk-\be_m}(\bx)\right.\\\notag
    &\left.+\sum_{l=1}^{d}4k_l(k_l-1)H_{\bk-2\be_l}(\bx)\right\} \\\notag    
    =&\sum_{\bk\in\Omega}\tilde{\varphi}_{j,\bk}^{i}\left\{-\sum_{l,m=1}^{d}\left(\Sigma_i^{-1}\right)_{lm}k_mH_{\bk-\be_m+\be_l}(\bx)\right.
       \end{align*}
    \begin{align} \label{eqn-3.15}	\notag
    &\left.+\sum_{l,m=1}^{d}\left(\Sigma_i^{-1}\right)_{lm}2k_mX_{t,l}^iH_{\bk-\be_m}(\bx)\right.\\\notag 
    &-\sum_{\substack{l,m=1\\l\not=m}}^{d}\left(\Sigma_i^{-1}\right)_{lm}2k_mk_lH_{\bk-\be_m-\be_l}(\bx)\\    
    &\left.+\sum_{l=1}^{d}2\left[2-\left(\Sigma_i^{-1}\right)_{ll}\right]k_l(k_l-1)H_{\bk-2\be_l}(\bx)\right\},
\end{align}
where the first and the third equalities are due to the properties of the Hermite polynomials, i.e. $H_k'(x)=2kH_{k-1}(x)$ and $H_{k+1}(x)=2xH_k(x)-2kH_{k-1}(x)$. Plugging \eqref{eqn-3.14} to \eqref{eqn-3.15}, multiplying both sides by $H_{\bq}(\bx)\bw(\bx)$, for $1\leq |\bq|_1\leq p$, and integrating with respect to $\bx\in\R^d$, it yields that
\begin{align}\label{eqn-3.16}\notag
	&C_j^i\bgamma_{\bf0}\bdelta_{{\bf0}\bq}-a_{j,\bq}\bgamma_{\bq}\\\notag
	\overset{\eqref{eqn-3.15}}=&-\sum_{l=1}^{d}\tilde{\varphi}_{j,\bq}^{i}\left(\Sigma_i^{-1}\right)_{ll}q_l\bgamma_{\bq}\\\notag
    &-\sum_{\substack{l,m=1\\\notag l\not=m}}^d\tilde{\varphi}_{j,\bq+\be_m-\be_l}^i\left(\Sigma_i^{-1}\right)_{lm}(q_m+1)\bgamma_{\bq}\\\notag
	&+\sum_{l,m=1}^d\tilde{\varphi}_{j,\bq+\be_m}^i\left(\Sigma_i^{-1}\right)_{lm}2(q_m+1)X_{t,l}^i\bgamma_\bq\\\notag
    &-\sum_{\substack{l,m=1\\\notag l\not=m}}^d\tilde{\varphi}_{j,\bq+\be_m+\be_l}^i\left(\Sigma_i^{-1}\right)_{lm}2(q_m+1)(q_l+1)\bgamma_\bq\\
	&+\sum_{l=1}^d\tilde{\varphi}_{j,\bq+2\be_l}^i2\left(2-\left(\Sigma_i^{-1}\right)_{ll}\right)(q_l+1)(q_l+2)\bgamma_{\bq},
\end{align}
since $\int_{\R^d}H_{\bk}(\bx)H_{\bq}(\bx)\bw(\bx)d\bx=\bgamma_\bq\bdelta_{\bk\bq}$, with $\bw(\bx):=\exp\left(-\sum_{l=1}^dX^2_l\right)$, $\bdelta_{\bk\bq}:=\prod_{l=1}^d\delta_{k_lq_l}$ and $\bgamma_{\bq}=\Pi_{l=1}^d\gamma_{q_l}$, $\gamma_{q_l}:=\pi^{\frac12}2^{q_l}q_l!$. Equation \eqref{eqn-3.17} is obtained by moving the terms with $\tilde\varphi_{j,\bq}^i$ and $\tilde\varphi_{j,\bq+\be_m-\be_l}^i$ to the left-hand side, and the rest to the right, for $1\leq|\bq|_1\leq p$. The constant $C_j^i$ is determined by substituting $\bq={\bf 0}$ in \eqref{eqn-3.16} and the convention that $\tilde\varphi_{j,\bq}^i\equiv0$, if some $q_l<0$, $l=1,\cdots,d$.
\end{IEEEproof}

\begin{remark}
	Let us write \eqref{eqn-3.17} as partitioned matrices:
\begin{align}\label{eqn-7}\notag
	\left(\begin{array}{cccccc}	B_1&D_2&&&&\\
	A_1&B_2&D_3&&&\\
        &&\ddots&\ddots&\ddots&\\
        &&&A_{p}&B_{p+1}&D_{p+2}
\end{array}\right)_{(p+1)\times(p+2)}\\
       \cdot \left(\begin{array}{c}
	\tilde\Phi_{j,1}^i\\
	\tilde\Phi_{j,2}^i\\
         \vdots\\
        \tilde\Phi_{j,p}^i\\
       \tilde\Phi_{j,p+1}^i\\
       \tilde\Phi_{j,p+2}^i\\
\end{array}\right)_{(p+2)\times 1}
=
\left(\begin{array}{c}
	{\bf a}_{j,0}-C_j^i\\
	{\bf a}_{j,1}\\
        \vdots\\
        {\bf a}_{j,p}
\end{array}\right)_{(p+1)\times1},
\end{align}
where $\tilde\Phi_{j,q}^i$ and ${\bf a}_{j,q}$ are the vectors of $\tilde\varphi_{j,\bq}^i$ and $a_{j,\bq}$, $|\bq|_1=q$, in some order, respectively,  $a_{j,\bq}$ are coefficients of $h_j$ in \eqref{eqn-3.14} , and $A_q$, $B_q$ and $D_q$ are the corresponding coefficient matrices in \eqref{eqn-3.17}, which depend on $\Sigma_i^{-1}$ and $\bq$s. The precise elements of $A_q$ are related to the order of $\tilde\varphi_{j,\bq}^i$, $|\bq|_1=q$, arranged in $\tilde\Phi_{j,q}^i$. Suppose that all $A_q$s are invertible, then starting from $\tilde\Phi_{j,p+1}^i={\bf 0}$ and $\tilde\Phi_{j,p+2}^i={\bf 0}$, backward recursively one has
\begin{align}\label{eqn-4}\notag
    \tilde\Phi_{j,p}^i=&A_{p}^{-1}{\bf a}_{j,p},\\\notag
    \tilde\Phi_{j,p-1}^i=&A_{p-1}^{-1}({\bf a}_{j,p-1}+B_{p}\tilde\Phi_{j,p}^i),\\\notag
	\tilde\Phi_{j,p-2}^i=&A_{p-2}^{-1}({\bf a}_{j,p-2}+B_{p-1}\tilde\Phi_{j,p-1}^i+D_{p}\tilde\Phi_{j,p}^i),\\
	\vdots&\\\notag
	\tilde\Phi_{j,1}^i=&A_1^{-1}({\bf a}_{j,1}+B_2\tilde\Phi_{j,2}^i+D_3\tilde\Phi_{j,3}^i).
\end{align}
At last, the constant $C_j^i$ is given by 
\[
    C_j^i={\bf a}_{j,{\bf 0}}-B_1\tilde\Phi_{j,1}^i-D_2\tilde\Phi_{j,2}^i.
\]
Therefore, the only question left during this procedure of backward recursion is the invertibility of $A_q$s, for all $1\leq q\leq p$. The invertibility will be discussed further in Section \ref{sec-3.2}. 
\end{remark}

We summarize the results in this subsection as following:

\begin{theorem}\label{thm-3.5}
	Under the assumption \eqref{eqn-3.14}, the gain function $\nabla\varphi_j(\bx)$ is obtained by 
\begin{align}\label{eqn-6}
	\nabla\varphi_j(\bx)
	=&\frac{1}{\sum_{i=1}^{N_p}\mathcal{N}(\bx;X_t^i,\Sigma_i)}\\\notag
    &\cdot\sum_{i=1}^{N_p}\left[(\bx-X_t^i)\frac{p_t^{N_p,\Sigma}[h_j]-C_j^i}{2\pi^{\frac{d}{2}}|\Sigma_i|^{\frac12}}\gamma\left(\frac d2,\frac{r^2}{2}\right)r^{-d}\right.\\\notag
    &\phantom{\sum_{i=1}^{N_p}aaz}\left.+\mathcal{N}(\bx;X_t^i,\Sigma_i)\sum_{\bk\in\Omega}\tilde\varphi_{j,\bk}^i2k_lH_{\bk-\be_l}(\bx)\right],
\end{align}
with the weighted radial $r:=\sqrt{\sum_{l=1}^{d}\lambda_l^i(x_l-X_{t,l}^i)^2}$, $\left\{\lambda_l^i\right\}_{l=1}^d$ being the eigenvalues of $\Sigma_i^{-1}$, such that $\nabla_{\lambda^i}\psi_j^i(X_t^i)={\bf 0}$, for $j=1,\cdots,m$, $i=1,\cdots,N_p$, $C_j^i$ is in \eqref{eqn-Cji}, and the coefficients $\tilde\varphi_{j,\bk}^i$, $1\leq|\bk|_1\leq p$, is obtained by the backward recursion \eqref{eqn-3.17}, with the initial condition $\tilde\varphi_{j,\bq}^i\equiv0$, for $|\bq|_1>p$, and the convention $\tilde\varphi_{j,\bq}^i\equiv0$, if some $q_l<0$. 
\end{theorem}

The next subsection is devoted to discuss the invertibility of the coefficient matrices $A_q$s in \eqref{eqn-7}. We believe that $A_q$s are invertible for all $d\geq1$ and $q\in\mathbb{N}$, but we can only prove this result in some special cases. The readers may skip Section \ref{sec-3.2} during their first reading, as this will not hinder their understanding of the content.

\subsection{The invertibility of the coefficient matrix in \eqref{eqn-7}}\label{sec-3.2}

In this subsection, we shall show the invertibility of the $A_q$s in some special cases.

\subsubsection{Case 1: for any $d,\ |\bq|_1=\sum_{l=1}^dq_l=q\in\mathbb{Z}_+$ \& the independent states}

\begin{proposition}\label{prop-3.6}
    For any $d,|\bq|_1=\sum_{l=1}^dq_l=q\in\mathbb{Z}_+$, and if the states are independent, i.e. $\Sigma_i=\textup{diag}(\lambda_1^i,\cdots,\lambda_d^i)$, for all $\lambda_l^i>0$, then the coefficient matrix $A_q$s in \eqref{eqn-7} are invertible, for $1\leq q\leq p$.
\end{proposition}
 \begin{IEEEproof}
    It is clear to see that $A_q$ is a diagonal matrix with the elements $\sum_{l=1}^d(\lambda_l^i)^{-1}q_l>0$, thus the invertibility follows immediately.
\end{IEEEproof}

\subsubsection{Case 2: for any $d\in\mathbb{Z}_+$, $|\bq|_1=1$ \& the dependent states}
\begin{proposition}\label{prop-3.7}
    For $d\in\mathbb{Z}_+$ arbitrary, $|\bq|_1=1$, and $\Sigma_i$ is arbitrary positive definite symmetric matrix, then $A_1$ in \eqref{eqn-7} are invertible.
\end{proposition}
 
 \begin{IEEEproof}
     All $\bq$s such that $|\bq|_1=1$ are $\be_k$, $k=1,\cdots,d$. If we arrange $\tilde{\varphi}_{j,\be_k}^i$ in the increasing order of $k$, it is easy to deduce that the matrix is exactly $\Sigma_i^{-1}$. In fact, if $\bq=\be_k$, for some $k=1,\cdots,d$, then the  two terms on the left-hand side of \eqref{eqn-3.17} are
\begin{align*}
\sum_{l=1}^d\tilde{\varphi}_{j,\be_l}^i\left(\Sigma_i^{-1}\right)_{ll}q_l=&\tilde{\varphi}_{j,\be_k}^i\left(\Sigma_i^{-1}\right)_{kk},
\end{align*}
and
\begin{align*}
	&\sum_{\substack{l,m=1\\\notag l\not=m}}^d\tilde{\varphi}_{j,\bq+\be_m-\be_l}^i\left(\Sigma_i^{-1}\right)_{lm}(q_m+1)\\\notag
	=&\sum_{\substack{m=1\\\notag m\not=k}}^d\tilde{\varphi}_{j,\be_m}^i\left(\Sigma_i^{-1}\right)_{km}(q_m+1)
	=\sum_{\substack{m=1\\\notag m\not=k}}^d\tilde{\varphi}_{j,\be_m}^i\left(\Sigma_i^{-1}\right)_{km},
\end{align*}
where the first equality follows by the fact that $l=k$, otherwise the $l$-th component of $\bq+\be_m-\be_l$ is negative, and the second equality is due to $q_m=0$, since $m\neq k$.
 \end{IEEEproof}
 
\subsubsection{Case 3: for  $d=2$, $|\bq|_1\in\mathbb{Z}_+$ \& the dependent states}  
\begin{proposition}\label{prop-3.9}
    For $d=2$ and $|\bq|_1\in\mathbb{Z}_+$, and $\Sigma_i$ is any positive definite symmetric matrix, then $A_p$s in \eqref{eqn-7} are invertible, $1\leq p\leq |\bq|_1$.
\end{proposition}

\begin{IEEEproof}
    The case $|\bq|_1=1$ has been discussed in Proposition \ref{prop-3.7}, case 2. Without loss of generality, we assume that $|\bq|_1\geq2$. All the possible $\bq$ such that $|\bq|_1=p\geq 2$ are $\{(p,0),(p-1,1),\cdots,(1,p-1),(0,p)\}$. If $\bq=p\be_k$, $k=1,2$, then the two terms on the left-hand side of \eqref{eqn-3.17} are
\begin{align*}
	\sum_{l=1}^2\tilde{\varphi}_{j,p\be_k}^i\left(\Sigma_i^{-1}\right)_{ll}q_l
	=&\tilde{\varphi}_{j,p\be_k}^ip\left(\Sigma_i^{-1}\right)_{kk},
\end{align*}
and
\begin{align*}
	&\sum_{\substack{m,l=1\\\notag l\not=m}}^2\tilde{\varphi}_{j,p\be_k+\be_m-\be_l}^i\left(\Sigma_i^{-1}\right)_{lm}(q_m+1)\\\notag
	=&\tilde{\varphi}_{j,(p-1)\be_k+\be_m}^i\left(\Sigma_i^{-1}\right)_{km},
\end{align*}
for $m\not=k$. If $\bq=p_1\be_1+p_2\be_2$, with $p_1,p_2\not=0$, $p_1+p_2=p$, then the two terms on the left-hand side of \eqref{eqn-3.17} becomes
\begin{align*}
	&\sum_{l=1}^2\tilde{\varphi}_{j,p_1\be_1+p_2\be_2}^i\left(\Sigma_i^{-1}\right)_{ll}q_l\\
    =&\tilde{\varphi}_{j,p_1\be_1+p_2\be_2}^i\left[\left(\Sigma_i^{-1}\right)_{11}p_1+\left(\Sigma_i^{-1}\right)_{22}p_2\right],
\end{align*}
and
\begin{align*}
	&\sum_{\substack{m,l=1\\\notag l\not=m}}^2\tilde{\varphi}_{j,p_1\be_1+p_2\be_2+\be_m-\be_l}^i\left(\Sigma_i^{-1}\right)_{lm}(q_m+1)\\
	=&\left(\Sigma_i^{-1}\right)_{12}(p_2+1)\tilde{\varphi}_{j,(p_1-1)\be_1+(p_2+1)\be_2}^i\\
    &+\left(\Sigma_i^{-1}\right)_{21}(p_1+1)\tilde{\varphi}_{j,(p_1+1)\be_1+(p_2-1)\be_2},
\end{align*}
since either $l=1,m=2$ or $l=2,m=1$. 

For the short of notation, let us denote the matrix $\Sigma_i^{-1}=\left(\begin{array}{cc}
	a&c\\
	c&b
\end{array}\right)$, with $a,b>0$ and \textup{det}$\left(\Sigma_i^{-1}\right)=ab-c^2>0$, due to the positive definiteness of $\Sigma_i$.   Consequently, the coefficient matrix $A_p$ infront of $\tilde\Phi_{j,p}^i=\left(\tilde\varphi_{j,(p,0)}^i,\cdots,\tilde\varphi_{j,(0,p)}^i\right)^T$ is tri-diagonal, i.e.
\begin{align}\label{eqn-a.2}
	\left(\begin{array}{ccccc}
	ap&c&0&&\\
	cp&a(p-1)+b&2c&&\\
	&c(p-1)&a(p-2)+2b&3c&\\
	&\ddots&\ddots&\ddots&\\
	&&0&c&bp
\end{array}\right)_{(p+1)\times(p+1)}.
\end{align}

We claim that \eqref{eqn-a.2} is invertible. The strategy is to perform elementary row operations to get an upper triangular matrix and to show that all the diagonal elements, denoted as $d_l$, $1\leq l\leq p+1$, are positive. One only need to show $d_l> a(p-l+1)+(l-1)\beta$, with $\beta:=b-\frac{c^2}a>0$, by mathematical induction, for $2\leq l\leq p+1$. Starting from $l=2$, we have 
$$d_2=a(p-1)+\beta>a(p-1)>0.$$ 
With the assumption that $d_l>a(p-l+1)+(l-1)\beta$, it yields that 
\begin{align*}
	d_{l+1}=&a(p-l)+lb+lc\left(-\frac{c(p-l+1)}{d_{l}}\right)\\
>&a(p-l)+lb+lc\left(-\frac{c(p-l+1)}{a(p-l+1)+(l-1)\beta}\right)\\
	>&a(p-l)+l\beta>0.
\end{align*}
\end{IEEEproof}

\subsubsection{Case 4: for $d=3$, $|\bq|_1=2$ \& the dependent states}
\begin{proposition}
    For $d=3$ and $|\bq|_1=2$, and $\Sigma_i$ is any positive definite symmetric matrix, then $A_1$ and $A_2$ in \eqref{eqn-7} are invertible.
\end{proposition}

\begin{IEEEproof}
    Notice that $A_1$ in \eqref{eqn-7} in front of $\tilde\Phi_{j,1}^i$ is invertible due to case 2. In the following, we only need to discuss the invertibility of $A_2$. When $|\bq|_1=2$, if $\bq=2\be_k$, $k=1,2,3$, then the two terms on the left-hand side of \eqref{eqn-7} are
\begin{align*}
\sum_{l=1}^3\tilde{\varphi}_{j,2\be_k}^i\left(\Sigma_i^{-1}\right)_{ll}q_l
=2\tilde{\varphi}_{j,2\be_k}^i\left(\Sigma_i^{-1}\right)_{kk},
\end{align*}
and
\begin{align*}
\sum_{\substack{m,l=1\\\notag l\not=m}}^3\tilde{\varphi}_{j,2\be_k+\be_m-\be_l}^i\left(\Sigma_i^{-1}\right)_{lm}(q_m+1)
=\sum_{\substack{l=1\\\notag l\not=k}}^3\tilde{\varphi}_{j,\be_k+\be_l}^i\left(\Sigma_i^{-1}\right)_{kl}.
\end{align*}
If $\bq=\be_l+\be_m$, with $l,m=1,2,3$, $l\not=m$, then the two terms on the left-hand side of \eqref{eqn-7} become
\begin{align*}
	\sum_{l=1}^3\tilde{\varphi}_{j,\be_l+\be_m}^i\left(\Sigma_i^{-1}\right)_{ll}q_l=\tilde{\varphi}_{j,\be_l+\be_m}^i\left[\left(\Sigma_i^{-1}\right)_{ll}+\left(\Sigma_i^{-1}\right)_{mm}\right],
\end{align*}
and
\begin{align*}
&\sum_{\substack{k,n=1\\\notag k\not=n}}^3\tilde{\varphi}_{j,\be_l+\be_m+\be_k-\be_n}^i\left(\Sigma_i^{-1}\right)_{lm}(q_m+1)\\
    =&2\left(\Sigma_i^{-1}\right)_{lm}\tilde{\varphi}_{j,2\be_m}^i+\left(\Sigma_i^{-1}\right)_{lk}\tilde{\varphi}_{j,\be_m+\be_k}
+2\left(\Sigma_i^{-1}\right)_{ml}\tilde{\varphi}_{j,2\be_l}^i\\
    &+\left(\Sigma_i^{-1}\right)_{mk}\tilde{\varphi}_{j,\be_l+\be_k},
\end{align*}
for $k\not=l$ or $m$. Let 
$\Sigma_i^{-1}\triangleq\left(\begin{array}{ccc}
a_{11}& a_{12}&a_{13}\\
a_{21}& a_{22}&a_{23}\\
a_{31}& a_{32}&a_{33}\\
\end{array}\right)$, then the coefficient matrix $A_2$ in \eqref{eqn-7} in front of $\tilde\Phi_{j,2}^i=[\tilde\varphi_{j,2\be_1}^i,\tilde\varphi_{j,2\be_2}^i,\tilde\varphi_{j,2\be_3}^i,\tilde\varphi_{j,\be_1+\be_2}^i,\tilde\varphi_{j,\be_1+\be_3}^i,\tilde\varphi_{j,\be_2+\be_3}^i]^T$ is
\begin{align}\label{eqn-11}\notag
	&A_2=\\
    &\left(\begin{array}{ccc:ccc}
	2a_{11} & 0 & 0 & a_{12} & a_{13} & 0\\
	0 & 2a_{22} & 0 & a_{12} & 0 & a_{23}\\
	0 & 0 & 2a_{33} & 0 & a_{13} & a_{23}\\
	\hdashline
	2a_{11} & 2a_{12} & 0 & a_{11}+a_{22} & a_{23} & a_{13}\\
	2a_{13} & 0 & 2a_{13} & a_{23} & a_{11}+a_{33} & a_{12}\\
	0 & 2a_{23} & 2a_{23} & a_{13} & a_{12} & a_{22}+a_{33}\\\notag
\end{array}\right)\\
&\triangleq\left(\begin{array}{c:c}
2A & B \\
\hdashline
2B^T & C \\ 
\end{array}\right).
\end{align}

To show the invertibility of $A_2$, it is equivalent to show $\textup{det}(A_2)>0$, i.e.
\begin{align}\label{}\notag
\textup{det}(A_2)&=\left|\begin{array}{cc}
2A & B \\
2B^T & C \\ 
\end{array}\right|
=2^3\textup{det}\left(AC-B^TB\right)\\
&=2^3\textup{det}(A)\cdot\textup{det}\left(C-B^TA^{-1}B\right).
\end{align}
$\textup{det}(A)=a_{11}a_{22}a_{33}>0$, since all the diagonal elements in the symmetric positive definite matrix are positive. Thus, we only need to prove $\textup{det}\left(C-B^TA^{-1}B\right)>0$. By direct computations, we have
\begin{align}\label{eqn-12}
	&C-B^TA^{-1}B 
	\overset{\eqref{eqn-11}}=\\\notag
    &\left(\begin{array}{ccc}
	\left(\frac{1}{a_{11}}+\frac{1}{a_{22}}\right)A_{33} & -\frac{1}{a_{11}}A_{23} & -\frac{1}{a_{22}}A_{13}\\
	-\frac{1}{a_{11}}A_{23} & \left(\frac{1}{a_{11}}+\frac{1}{a_{33}}\right)A_{22} & -\frac{1}{a_{33}}A_{12} \\
	-\frac{1}{a_{22}}A_{13} & -\frac{1}{a_{33}}A_{12}  & \left(\frac{1}{a_{22}}+\frac{1}{a_{33}}\right)A_{11}
\end{array}\right),
\end{align}
where $A_{lk}$ is the $(l,k)$-th algebraic cofactor of $\Sigma_i^{-1}$, $l,k=1,2,3$. Let us rewrite 
\begin{equation}\label{eqn-20}
	C-B^TA^{-1}B:=\bar A_1+\bar A_2+\bar A_3,
\end{equation}
 where
\begin{align*}\label{}
\bar{A}_1:=\frac{1}{a_{11}}\left(\begin{array}{ccc}
A_{33} & -A_{23} & 0\\
-A_{23} & A_{22} & 0 \\
0 & 0 & 0 \\
\end{array}\right),
\end{align*}
\begin{align*}
\bar{A}_2:=\frac{1}{a_{22}}\left(\begin{array}{ccc}
A_{33} & 0 & -A_{13}\\
0 & 0 & 0 \\
-A_{13} & 0 & A_{33} \\
\end{array}\right),
\end{align*}
\begin{align*}
	\bar{A}_3:=\frac{1}{a_{33}}\left(\begin{array}{ccc}
0 & 0 & 0\\
0 & A_{22} & -A_{12} \\
0 & -A_{12} &A_{33} \\
\end{array}\right).
\end{align*}
We claim that the submatrices of $\bar A_i$ with zero elements removed, denoted as $\bar A_i'$ are all positive definite, $i=1,2,3$. 

Let us show the positive definiteness of $\bar A_1'=\frac{1}{a_{11}}\left(\begin{array}{cc}
A_{33} & -A_{23}\\
-A_{23} & A_{22}
\end{array}\right)$ as an example. Notice that given that $\Sigma_i^{-1}$ is positive definite, so is the adjoint matrix $\left(\Sigma_i^{-1}\right)^{*}=\left(\begin{array}{ccc}
A_{11}& A_{21} & A_{31}\\
A_{12} & A_{22} & A_{32} \\
A_{13} & A_{23} &A_{33} \\
\end{array}\right)$. All its principal minors are positive, for example $A_{33}>0$, $\left|\begin{array}{cc}
A_{22} & A_{32} \\
A_{23} &A_{33} \\
\end{array}\right|=A_{22}A_{33}-A_{23}A_{32}>0$. Thus, all the sequential principal minors of $\bar A_1'$ are positive, since $\left|\begin{array}{cc}
A_{33} & -A_{23} \\
-A_{32} &A_{22} \\
\end{array}\right|=A_{22}A_{33}-A_{23}A_{32}>0$. Consequently, $\bar A_1'$ is positive definite. The same arguments yield that $\bar A_2'$ and $\bar A_3'$ are also positive definite. 

By the definition of positive definite, for all $\boldsymbol{v}\not={\bf 0}\in\R^2$, $\boldsymbol{v}\bar{A}_i'\boldsymbol{v}^T>0$, $i=1,2,3$. Therefore, for all $\boldsymbol{v}\not={\bf 0}\in\R^3$, there exists at least one $\bar A_i$ such that $\boldsymbol{v}\bar A_i\boldsymbol{v}^T>0$. Therefore, for all $\boldsymbol{v}\not={\bf 0}\in\R^3$, $\boldsymbol{v}\left(C-B^TA^{-1}B\right)\boldsymbol{v}^T\overset{\eqref{eqn-20}}=\boldsymbol{v}\left(\bar A_1+\bar A_2+\bar A_3\right)\boldsymbol{v}^T>0$, which asserts that $C-B^TA^{-1}B$ is positive definite.
\end{IEEEproof}

For more general cases, we conjecture that the matrices $A_q$s
in \eqref{eqn-7} are invertible. While we intend to further address this issue in future work, all the special cases discussed in this subsection already cover the examples presented in Section \ref{sec-4}.

\section{Numerical simulations}\label{sec-4}

\subsection{Algorithm}\label{AF}

\IEEEPARstart{T}{o} verify the effectiveness of our decomposition algorithm, we employ Stratonovich-form filters \cite{YMMS:12} and Euler discretization. The discrete-time algorithm derived therefrom is presented in Algorithm \ref{alg-1}.

\begin{algorithm}[ht!]
	\caption{The decomposition method for FPF}
	\label{alg-1}
	\begin{algorithmic}[1]
		\STATE{\textbf{\%Initialization}}
		\FOR{$i=1$ to $N_p$}
		\STATE{Sample $X_0^i$ from $p_0(x)$}
		\ENDFOR
		\STATE{\textbf{\% The FPF}}
        \FOR{$k=0$ to $k=T/\Delta t$}
		\STATE Calculate $p_{t_{k}}^{N_p,\Sigma}[h_j]\approx\frac{1}{N_p}\sum_{i=1}^{N_p}h_j\left(X_{t_k}^{i}\right)$, $j=1,\cdots,m$.
		\FOR{$i=1$ to $N_p$} 
            \STATE Generate $N_p$ independent samples $\Delta B^{i}_{t_k}$ and $\Delta W_{t_k}^i$ from $\N(0,\Delta t)$, respectively
            \STATE {\textbf{\% The decomposition method}}
            \STATE  The coefficients $\tilde\varphi_{j,\bk}^i$, $1\leq|\bk|_1\leq p$, is obtained by the backward recursion \eqref{eqn-3.17}, and $C_j^i$ is in \eqref{eqn-Cji}
            \STATE           
            Calculate $K_{\cdot j}=\nabla\varphi_j$ and $u$ for the $i$-th particle $X_{t_k}^{i}$ by \eqref{eqn-6} and \eqref{eqn-u}, respectively
		  \STATE Evolve the particles 
            $\{X_{t_k}^{i}\}_{i = 1}^{N_p}$ according to \eqref{eqn-2.20}-\eqref{eqn-2.3}, i.e.
		  \begin{align*}
            X_{t_{k+1}}^i=&X_{t_k}^i+g\left(X_{t_k}^i\right)\Delta t+\sigma\left(X_{t_k}^i\right)\sqrt{\Delta t}\Delta B_{t_k}^i\\
            &+u\left(X_{t_k}^i,t_k\right)\Delta t+K\left(X_{t_k}^i,t_k\right)\Delta Z_{t_k}^i,
            \end{align*}
            where $\Delta Z_{t_k}^i=h(X_{t_k}^i)\Delta t+\Delta W_{t_k}^i$.
		\ENDFOR
        \ENDFOR
	\end{algorithmic}
\end{algorithm}


\subsection{The comparison of the gain functions}

We shall first compare the gain functions obtained by various methods with the exact one in a simple multivariate example.

\begin{example}\label{ex-4.1}
Suppose the true distribution $p_t^*$, at some time instant $t$, is a mixture of two Gaussian distributions in $\R^2$, given by $\frac12\N(\bx;\vec{\mu},\Sigma)+\frac12\N(\bx;-\vec{\mu},\Sigma)$, where $\vec{\mu}=(\mu_1,\mu_2)^T$ and $\Sigma=\textup{diag}(\sigma_1,\sigma_2)$. The observation is a scalar-valued function: $h(\bx)=\frac{1}{\sigma_1}x_1^2-\frac{\mu_1}{\sigma_1\mu_2}x_1x_2$. Through a direct computation, one can easily verify that the gain function is
$K(\bx)=\left(x_1,-\frac{\sigma_2\mu_1}{\sigma_1\mu_2}x_1\right)^T$.
 
\begin{table}[bh!]
    \centering
\begin{tabular}{|c|c|c|c|}
\hline
          $l^2(K_{\cdot,j},K_j)$ & $\cdot=const$ & $\cdot=ker$ &$\cdot=decomp$\\
    \hline
     $j=1$ & $16.1567$& $22.9869$&$\bf{9.1351}$ \\
     \hline
     $j=2$ & $31.7946$ & $34.0134$ &$\bf{9.5147}$ \\
    \hline
\end{tabular}
\vspace{1em}
    \caption{The $l^2$-errors of the gain functions obtained by three different methods}
    \label{table-0}
\end{table}
We shall obtain the approximate gain functions by three different methods:  the constant-gain approximation, the kernel-based method \cite{TM:16} and the decomposition method, denoted as $K_{const}$, $K_{ker}$ and $K_{decomp}$, respectively. The parameters in the decomposition method are set to be $N_p=100$, $\mu_1=\mu_2=1$, $\Sigma=\textup{diag}(\sigma_1,\sigma_2)$ with $\sigma_1=1$, $\sigma_2=2$ and $\Sigma_i=\textup{diag}(\varepsilon_1,\varepsilon_2
)$ with $\varepsilon_1=0.5$, $\varepsilon_2=1$. In the kernel-based method, the bandwidth of kernel is set to be $\varepsilon = 0.1$ and the number of iterations to approximate the semi-group operator is $T_{iter} = 100$. The performances are measured by the $l^2$-error of $K$ component-wisely, defined as
\begin{equation*}    l^2(K_{\cdot,j},K_j)=\sqrt{\sum_{i=1}^{N_p}\left[K_{\cdot,j}(X^i_t)-K_j(X^i_t)\right]^2},
\end{equation*}
$j=1,\cdots,d$, where $\{X_t^i\}_{i=1}^{N_p}$ are independent particles drawn from $p_t^*$, $K_{\cdot,j}$ represents the $j$th component of the approximate gain function $K$ obtained by one of the three methods $\cdot$ and $K_j$ is the $j$th component of the true one. The $l^2$-errors are list in Table \ref{table-0}, which reveals a significant improvement in both $K_1$ and $K_2$ by using the decomposition method. Notably, $K_1$ exhibits higher accuracy than $K_2$ across all three methods, though this difference is less pronounced in the decomposition method. This observation can be attributed to the true gain function in this example depending solely on the first variable $x_1$, rather than $x_2$. 

To provide a clearer visualization, we plot the approximate gain functions for both components against the variables separately in Fig. \ref{fig-1}. Each dot in the figure represents the gain function's value at each particle $X^i=(X_1^i,X_2^i)$, $i=1,\cdots,100$. Both $K_1$ and $K_2$ exhibit an almost linear relationship with $x_1$, whereas they fluctuate around a near-zero constant when plotted against $x_2$. Among the three methods, the decomposition method yields the closest approximation, while the kernel-based method generates more outliers.

\begin{figure}[ht!]
\setlength{\abovecaptionskip}{4pt}
\centering
\begin{minipage}[t]{0.5\textwidth}
    \centering 
    \includegraphics[width=\textwidth, trim = 10 200 10 180,clip]{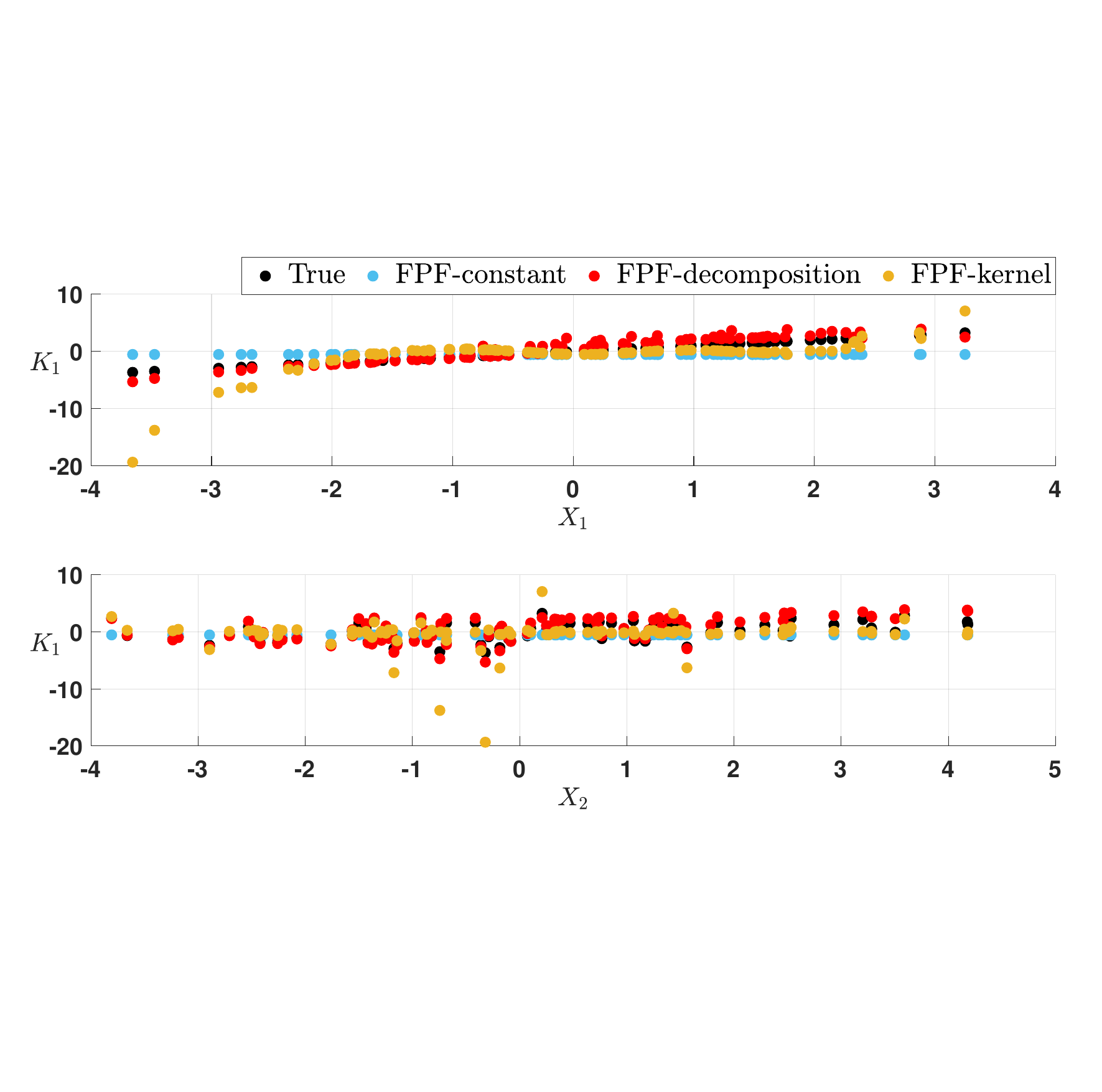}
    \small 
\end{minipage}
\vspace{4pt}
\begin{minipage}[t]{0.5\textwidth}
    \centering  
    \includegraphics[width=\linewidth, trim = 10 200 10 220,clip]{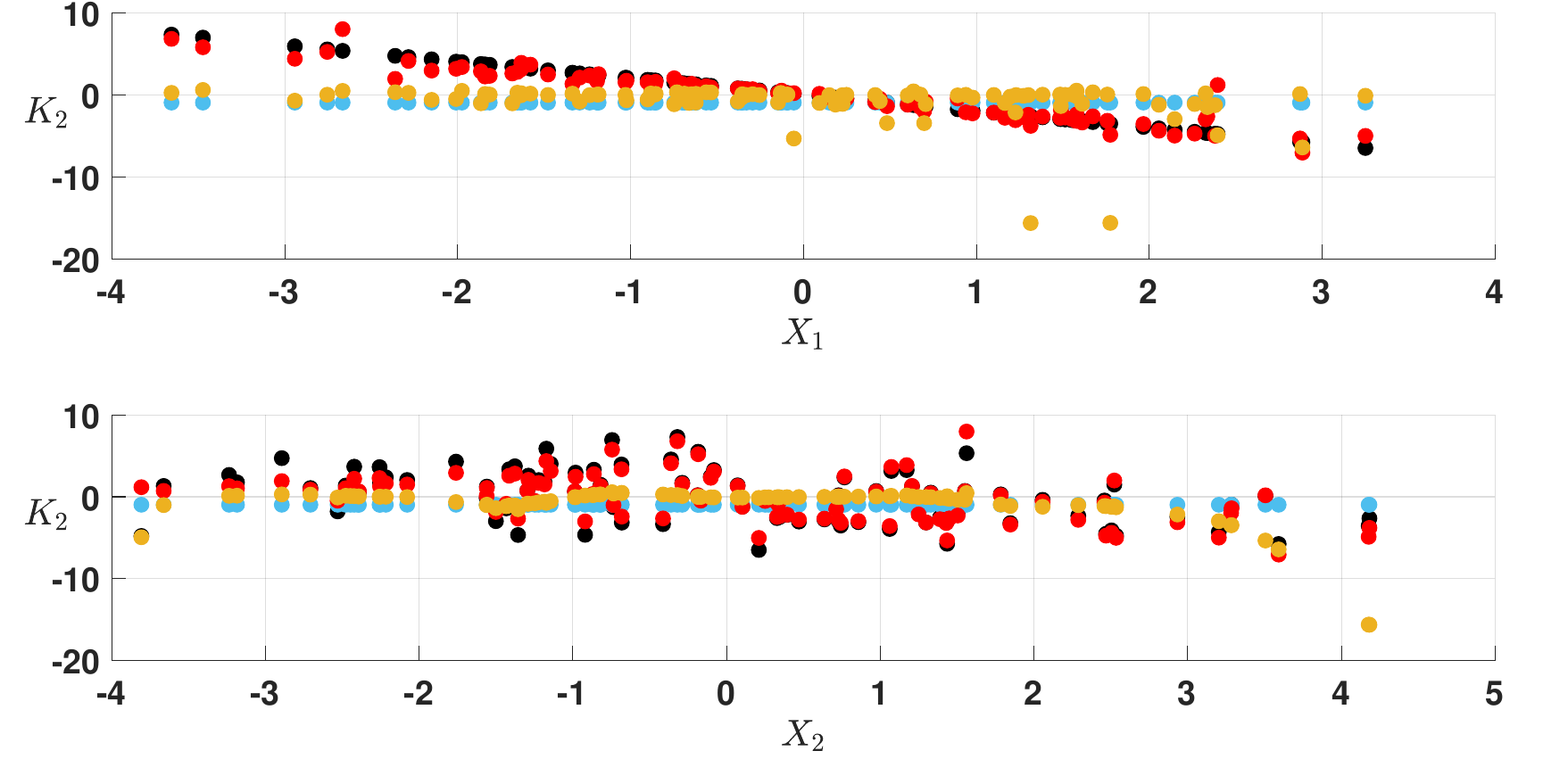}
    \small 
\end{minipage}
\small
\caption{Comparisons of the true gain function with those obtained by the three methods}
\label{fig-1}
\end{figure}
\end{example}

Next, we investigate the computational complexity of the multivariate decomposition method with respect to the state dimension. We design an essentially $d$-dimensional decoupled NLF problem. By virtue of its decoupling property, the gain function for each component $X_j$ can be readily obtained via Theorem III.5 in \cite{WML:25}. Meanwhile, the gain function can also be derived directly from Theorem \ref{thm-3.5} without relying on this decoupling property. We first derive these two distinct gain functions-they are not identical, yet both yield the same feedback to each particle $X_t^i$. This is because each particle is controlled by the gain function’s value evaluated at that particle. We then discuss and assess the efficiency and accuracy of the multivariate decomposition method in high-dimensional scenarios.

\begin{example}\label{ex-4.2}
Let us consider an extension of the cubic sensor problem to the multivariate dimensional situation:
\begin{equation}\label{eqn-4.7}
\left\{\begin{aligned}
dX_t=&X_t\odot(1_{d\times1}-X_t\odot X_t)dt+dB_t\\
dZ_t=&h(X_t)dt+dW_t
\end{aligned}\right.,
\end{equation}
where $1_{d\times1}$ is a $d \times 1$ column vector with all $1$s and $\odot$ is the Hadamard product of two vectors. The observation function is set to be $h(\bx)=(h_1(\bx),\cdots,h_m(\bx))^T$, with $m=d$ and $h_j(\bx)=x_j^3$, $j=1,\cdots,d$, $x_j$ is the $j$th component of $\bx$.  

Notice that \eqref{eqn-4.7} in fact is a $d$ decoupled scalar NLF problems. We have already investigated the decomposition method for the gain function in the scalar NLF problems in \cite{WML:25}. Here, we take advantage of this special setting to get the gain function by Theorem III.5, \cite{WML:25} and compare it with that obtained by the decomposition method for multivariate FPF, Theorem \ref{thm-3.5}. In particular, for each observation function 
\begin{equation}\label{eqn-4.9}
    h_j(\bx)=x_j^3=\frac1{8}H_3(x_j)+\frac3{4}H_1(x_j),
    \end{equation}
$j=1,\cdots,d$, where $H_q(\cdot)$ is the scalar Hermite polynomials of degree $q$. The backward recursion starting from $q=4$ is performed down to $q=0$, with the initial conditions $\tilde K_{j,4}^i=\tilde K_{j,3}^i=0$. By Theorem III.5, \cite{WML:25}, the backward recursion yields that
\begin{align}\label{eqn-4.18}
    \tilde K^i_{j,2}=&\frac{\varepsilon}{4},\quad
    \tilde K^i_{j,1}=\frac{\varepsilon X_{t,j}^i}{2},\quad
    \tilde K^i_{j,0}=\frac{\varepsilon}{2}+2\varepsilon^2+ {X_{t,j}^i}^2\varepsilon,
\end{align}
and $C_j^i={X_{t,j}^i}^3+3X_{t,j}^i\varepsilon$. In the notations of $\tilde\varphi_{j,q}^i$ in \eqref{eqn-3.4}, one has the relation
\begin{equation}\label{eqn-4.17}
    2(q+1)\tilde\varphi_{j,q+1}^i=\tilde K_{j,q}^i,
\end{equation}
$q=0,1,\cdots,p-1$, $j=1,\cdots,d$ and $i=1,\cdots,N_p$. Consequently, by Theorem III.5 in \cite{WML:25}, the gain function for the multivariate FPF is $K(\bx)=(K_1(x_1),\cdots,K_d(x_d))$, where $K_j(x_j)$ has the explicit expression 
\begin{align}\label{eqn-4.8}\notag
&K_j(x_j)
=\frac{1}{\sum\limits_{i=1}^{N_p}\mathcal{N}(x_j;X_{t,j}^i,\varepsilon)}\\
&\cdot\sum\limits_{i=1}^{N_p}\left\{\mathcal{N}(x_j;X_{t,j}^i,\varepsilon)(\varepsilon x_j^2+X_{t,j}^i\varepsilon x_j+2\varepsilon^2+{X_{t,j}^i}^2\varepsilon)\right.\\\notag
&\phantom{aaaa}+\left.\frac12\left[p_t^{N_p,\varepsilon}[h_j]-\left({X_{t,j}^i}^3+3X_{t,j}^i\varepsilon\right)\right]\textup{erf}\left(\frac{x_j-X_{t,j}^i}{\sqrt{2\varepsilon}}\right)\right\},
\end{align}
where $p_t^{N_p,\varepsilon}[h_j]:=\int_{\R}h_j(x_j)p_t^{N_p,\varepsilon}(x_j)dx_j$.

In the below, we detail the gain function $K_j(\bx)=\nabla\varphi_j(\bx)$, $j=1,\cdots,d$, by Proposition \ref{P2} and Theorem \ref{thm-3.5} without taking advantage of the decoupled-ness. According to \eqref{eqn-4.9}, given $j$ in Assumption \eqref{eqn-3.14} the coefficients $a_{j,\bq}\equiv0$, except $a_{j,3\be_j}=\frac18$ and $a_{j,\be_j}=\frac34$. By Remark \ref{rmk-10}, setting $p=3$, $\tilde\Phi_{j,4}^i={\bf0}$ and $\tilde\Phi_{j,5}^i={\bf0}$, the coefficients $\tilde\Phi_{j,q}^i$ satisfy:
\begin{align*}
    	&\left(\begin{array}{cccccc}	B_1&D_2&&&&\\
	A_1&B_2&D_3&&&\\
        &&A_2&B_3&D_4&\\
        &&&A_3&B_4&D_5
\end{array}\right)_{4\times5}\,\left(\begin{array}{c}
	\tilde\Phi_{j,1}^i\\
	\tilde\Phi_{j,2}^i\\
       \tilde\Phi_{j,3}^i\\
      {\bf 0}\\
       {\bf 0}\\
\end{array}\right)_{5\times 1}\\
&\phantom{aaaaaaaa}=
\left(\begin{array}{c}
	-C_j^i\\
	{\bf a}_{j,1}\\
    0\\
    {\bf a}_{j,3}
\end{array}\right)_{4\times1},
\end{align*}
where ${\bf0}$ is a zero vector of proper dimension, ${\bf a}_{j,q}$ is the vector of $a_{j,\bq}$, $|\bq|_1=q$, in some order, so do $\tilde\Phi_{j,q}$s. Clearly, the backward recursion is given as follows:
\begin{align}\label{eqn-4.10}
    \tilde\Phi_{j,3}^i=&A_{3}^{-1}{\bf a}_{j,3},\\\label{eqn-4.11}
    \tilde\Phi_{j,2}^i=&-A_{2}^{-1}B_{3}\tilde\Phi_{j,3}^i,\\\label{eqn-4.12}
	\tilde\Phi_{j,1}^i=&-A_{1}^{-1}(-{\bf a}_{j,1}+B_{2}\tilde\Phi_{j,2}^i+D_{3}\tilde\Phi_{j,3}^i)
\end{align}
and 
\begin{equation}\label{eqn-4.13}
    C_j^i=-(B_1\tilde\Phi_{j,1}^i+D_2\tilde\Phi_{j,2}^i),
\end{equation}
where $A_q$, $B_q$ and $D_q$ are closely related to the order of the index $\bq$, $|\bq|_1=0,1,2,3$. Starting from \eqref{eqn-4.10}, we find the coefficient $\tilde\Phi_{j,q}^i$ one by one:
\begin{enumerate}
    \item For $\tilde\Phi_{j,3}^i$: By Proposition \ref{prop-3.6}, $A_3=\textup{diag}\left(\frac3{\varepsilon} 1_{\frac{d(d+1)(d+2)}6\times1}\right)$, since $\sharp\{\bq:\,|\bq|_1=3\}=d+2 {d\choose 2}+{d\choose 3}=\frac{d(d+1)(d+2)}6$. Notice that ${\bf a}_{j,3}$ with all zero components except $a_{j,3\be_j}=\frac18$. Let us arrange $\bq$s starting with $(3\be_1,\cdots,3\be_d)$ and following with all the remaining $2\be_{l}+\be_m$, $\be_l+\be_m+\be_n$, $l\neq m\neq n$ in arbitrary order. Thus, \eqref{eqn-4.10} can be written as
  \begin{align}\label{eqn-4.20}\notag
	&  \left(\begin{array}{c:c}
    \frac3\varepsilon\textup{diag}(1_{d\times d}) & {\bf 0}_{d\times\frac{d(d-1)(d+4)}6}\\
	\hdashline
	  {\bf 0}_{\frac{d(d-1)(d+4)}6\times d}&\frac3\varepsilon\textup{diag}(1_{\frac{d(d-1)(d+4)}6\times\frac{d(d-1)(d+4)}6})
      \end{array}\right)\\
        &\phantom{aaaaaa}\cdot\left(\begin{array}{c}
    \tilde\varphi_{j,3\be_1}^i\\
    \vdots\\
    \tilde\varphi_{j,3\be_j}^i\\
    \vdots\\
    \tilde\varphi_{j,3\be_d}^i\\
\hdashline
    \textup{other}\\
    \tilde\varphi_{j,\bq}^i s
\end{array}\right)
    =\left(\begin{array}{c}
    0\\
    \vdots\\
    \frac18\\
    \vdots\\
    0\\
\hdashline
    {\bf 0}_{\frac{d(d-1)(d+4)}6\times1}
\end{array}\right).
\end{align}
This yields that all $\tilde\varphi_{j,\bq}^i=0$, with $|\bq|_1=3$, except 
\begin{equation}\label{eqn-4.14}
    \tilde\varphi_{j,3\be_j}^i\overset{\eqref{eqn-4.10}}=\frac\varepsilon{24}.
\end{equation}
    \item For $\tilde\varPhi_{j,2}^i$: again by Proposition \ref{prop-3.6}, $A_2=\textup{diag}\left(\frac2{\varepsilon} 1_{\frac{d(d+1)}2}\right)$, since $\sharp\{\bq:\,|\bq|_1=2\}={d\choose 2}+d=\frac{d(d+1)}2$. Notice from 1) that only $\tilde\varphi_{j,3\be_j}^i\neq0$, according to \eqref{eqn-4.11}, it is sufficient to find the corresponding component in $B_3$. For all $\bq$, $|\bq|_1=2$, the term
    \[
    -\sum_{l,m=1}^d\frac1\varepsilon2(q_m+1)X_{t,l}^i\tilde\varphi_{j,\bq+\be_m}^i=-\frac6\varepsilon X_{t,j}^i\tilde\varphi_{j,3\be_j}^i\neq0,
    \]
    only if $\bq=2\be_j$ and $l=m=j$. Consequently, all $\tilde\varphi_{j,\bq}^i=0$, with $|\bq|_1=2$, except
    \begin{equation}\label{eqn-4.15}
        \tilde\varphi_{j,2\be_j}^i\overset{\eqref{eqn-4.11},\eqref{eqn-4.14}}=-\frac\varepsilon2 \left(-\frac{6X_{t,j}^i}\varepsilon\right)\frac\varepsilon{24}=\frac{\varepsilon X_{t,j}^i}8.
    \end{equation}
    \item For $\tilde\varphi_{j,1}^i$: again by Proposition \ref{prop-3.6}, $A_1=\textup{diag}\left(\frac1\varepsilon1_{d\times1}\right)$, since $\sharp\{\bq:\,|\bq|_1=d\}$. Same as in 2), notice that $\tilde\varphi_{j,\bq}^i=0$ except $\tilde\varphi_{j,3\be_j}^i\neq0$. According to \eqref{eqn-4.12}, it is sufficient to find the corresponding component in $D_3$. For all $\bq$ with $|\bq|_1=1$, the term
    \begin{align*}
        &-\sum_{\substack{l,m=1\\l\not=m}}^d\tilde{\varphi}_{j,\bq+\be_m+\be_l}^i\left(\Sigma_i^{-1}\right)_{lm}2(q_m+1)(q_l+1)\\
	&+\sum_{l=1}^d\tilde{\varphi}_{j,\bq+2\be_l}^i2\left[2-\left(\Sigma_i^{-1}\right)_{ll}\right](q_l+1)(q_l+2)\\
    =&-12\left(2-\frac1\varepsilon\right)\tilde\varphi_{j,3\be_j}^i,
    \end{align*}
    only if $\bq=\be_j$ and $l=j$, where the first term vanishes due to the fact that $\bq+\be_m+\be_l\neq3\be_j$, if $|\bq|_1=1$ and $m\neq l$. Similarly, notice that only $\tilde\varphi_{j,2\be_j}^i\neq0$. According to \eqref{eqn-4.12}, it is sufficient to find the corresponding component in $B_2$. For all $\bq$ with $|\bq|_1=1$, the term
     \[
    -\sum_{l,m=1}^d\frac1\varepsilon2(q_m+1)X_{t,l}^i\tilde\varphi_{j,\bq+\be_m}^i=-\tilde\varphi_{j,2\be_j}^i\frac{4X_{t,j}^i}\varepsilon\neq0,
    \]
    only if $\bq=\be_j$ and $l=m=j$. Therefore, we have $\tilde\varphi_{j,\bq}^i=0$, $|\bq|_1=1$, except
    \begin{align}\notag
        \tilde\varphi_{j,\be_j}^i\overset{\eqref{eqn-4.12}}=&-\varepsilon\left(-\frac{4X_{t,j}^i}\varepsilon\frac{\varepsilon X_{t,j}^i}8-12\left(2-\frac1\varepsilon\right)\frac\varepsilon{24}-\frac34\right)\\\label{eqn-4.16}
        \overset{\eqref{eqn-4.14},\eqref{eqn-4.15}}=&\frac{\varepsilon {X_{t,j}^i}^2}2+\varepsilon^2+\frac\varepsilon4.
    \end{align}
    \item For $C_j^i$: Notice that only $\tilde\varphi_{j,2\be_j}^i\neq0$ and $\tilde\varphi_{j,\be_j}^i\neq0$, from \eqref{eqn-4.13}, it is sufficient to determine the corresponding component in $D_2$ and $B_1$, respectively. Similarly, as in 3), one has
    \begin{align*}
        C_j^i\overset{\eqref{eqn-4.13}}=&-\left[\left(-\frac{2X_{t,j}^i}\varepsilon\right)\tilde\varphi_{j,2\be_j}^i+\left(-4\left(2-\frac1\varepsilon\right)\tilde\varphi_{j,\be_j}^i\right)\right]\\
        \overset{\eqref{eqn-4.15},\eqref{eqn-4.16}}=&{X_{t,j}^i}^3+3\varepsilon X_{t,j}^i.
    \end{align*}
\end{enumerate}

Recall the relation between $\tilde K_{j}^i$ and $\tilde\varphi_{j,\bq}^i$ \eqref{eqn-4.17}, \eqref{eqn-4.14}-\eqref{eqn-4.16} match exactly with those in \eqref{eqn-4.18}, so does $C_j^i$. By Theorem \ref{thm-3.5}, the gain function $K(\bx)=(\nabla\varphi_1(\bx),\cdots,\nabla\varphi_d(\bx))$, where
\begin{align}\label{eqn-4.19}
    \nabla&\varphi_j(\bx)
	=\frac{1}{\sum_{i=1}^{N_p}\mathcal{N}(\bx;X_t^i,\varepsilon I_d)}\\\notag
    &\cdot\sum_{i=1}^{N_p}\left[\mathcal{N}(\bx;X_t^i,\varepsilon I_d)(\varepsilon \bx_j^2+X_{t,j}^i\varepsilon \bx_j+2\varepsilon^2+{X_{t,j}^i}^2\varepsilon)\right.\\\notag
    &\phantom{\sum_{i=1}^{N_p}aaz}\left.+(\bx-X_t^i)\frac{p_t^{N_p,\varepsilon I_d}[h_j]-C_j^i}{2(\pi\varepsilon)^{\frac{d}{2}}}\gamma\left(\frac d2,\frac{r^2}{2}\right)r^{-d}\right],
\end{align}
where $r=\sqrt{\frac1\varepsilon\sum_{l=1}^d(x_l-X_{t,l}^i)^2}$. Comparing the gain functions in \eqref{eqn-4.8} and \eqref{eqn-4.19}, $K_j\neq\nabla\varphi_j$, for $\bx\in\R^d$. Nevertheless, they are the same at $\bx=X_t^i$, since both $\tilde K_0(\bx)$ in \eqref{eqn-2.2} and $\nabla_{\lambda^i}\psi_j^i$ in \eqref{eqn-15} vanish. This makes the feedback to every particles being the same.

Take a closer look, the essential reason for this difference is due to the fact that not the Gaussian mixture $p_t^{Np,\varepsilon I_d}$ \eqref{eqn-2.4} is used to get \eqref{eqn-4.8}. In fact, in the decoupled derivation, the proposed approximation for each component $x_j$ is a one-dimensional Gaussian mixture, i.e.
\[
    p_t^{N_p,\varepsilon}(x_j)=\frac1{N_p}\sum_{i=1}^{N_p}\N(x_j;X_{t,j}^i,\varepsilon),
\]
for $j=1,\cdots,d$, with the assumption that the components are mutually independent, then the joint distribution is
\begin{align}\label{eqn-4.21}
    \tilde p_t^{N_p,\varepsilon I_d}(\bx)=&\prod_{j=1}^d p_t^{N_p,\varepsilon}(x_j)
        =\prod_{j=1}^d\frac1{N_p}\sum_{i=1}^{N_p}\N(x_j;X_{t,j}^i,\varepsilon),
\end{align}
which is different from the multivariate Gaussian mixture \eqref{eqn-4}:
\begin{equation}\label{eqn-4.22}
    p_t^{N_p,\varepsilon I_d}(\bx)=\frac1{N_p}\sum_{i=1}^{N_p}\prod_{j=1}^d\N(\bx_j;X_{t,j}^i,\varepsilon).
\end{equation}

Obviously, \eqref{eqn-4.21} and \eqref{eqn-4.22} are the same when $d=1$ or $N_p=1$. In general, the multivariate Gaussian mixture in \eqref{eqn-4} is more suitable for modeling when no prior knowledge of independence exists, as the covariance matrix $\Sigma$ can readily characterize the correlation among components $x_j$. However, if decoupling is known a priori, it should be fully exploited, as it can dramatically reduce computational burden-down to $\mathcal{O}(pd)$, as illustrated in Example IV.2 of \cite{WML:25}.

We analyze the computational complexity of \eqref{eqn-4.19} when the coefficients $\tilde{\Phi}_{j,q}^i$, $q=1,\cdots,p$ are non-sparse. It is straightforward to observe that inverting the matrices $A_q$, $q=1,\cdots,p$, constitutes the most computationally intensive step. Even for $\Sigma=\varepsilon I_d$, this inversion requires $\mathcal{O}(pd^3)$ operations per $j=1,\cdots,d$ by \eqref{eqn-4.20}, leading to a theoretical computational complexity of $\mathcal{O}(pd^4)$ for the gain function in \eqref{eqn-4.19}. For dense $A_q$ matrices (e.g., when $\Sigma_i$ are non-diagonal), the complexity increases to $\mathcal{O}(pd^7)$. Nevertheless, this polynomial scaling represents a significant improvement over the exponential growth of most particle-based algorithms. We conduct numerical experiments to investigate the relationship between CPU time and state dimension $d$ for system \eqref{eqn-4.7} using \eqref{eqn-4.19}, with parameters $T=40$, $\Delta t=0.01$, $N_p=50$ particles, $\Sigma_i=\varepsilon I_d$, $\varepsilon=0.01$, and $d=1,3,5,10,20,30,50,70,100$. As shown in Fig. \ref{cpu_time}, the log-transformed CPU time exhibits a polynomial relationship with the log-transformed dimension $d$, with a degree of $2.88$. We attribute this complexity (slightly lower than the theoretical degree $4$) to the sparsity of the coefficients $\tilde{\Phi}_{j,q}^i$ in the specific example \eqref{eqn-4.7}.
\begin{figure}[th!]
	\centering
        \includegraphics[width=0.5\textwidth, trim = 20 180 20 220,clip]{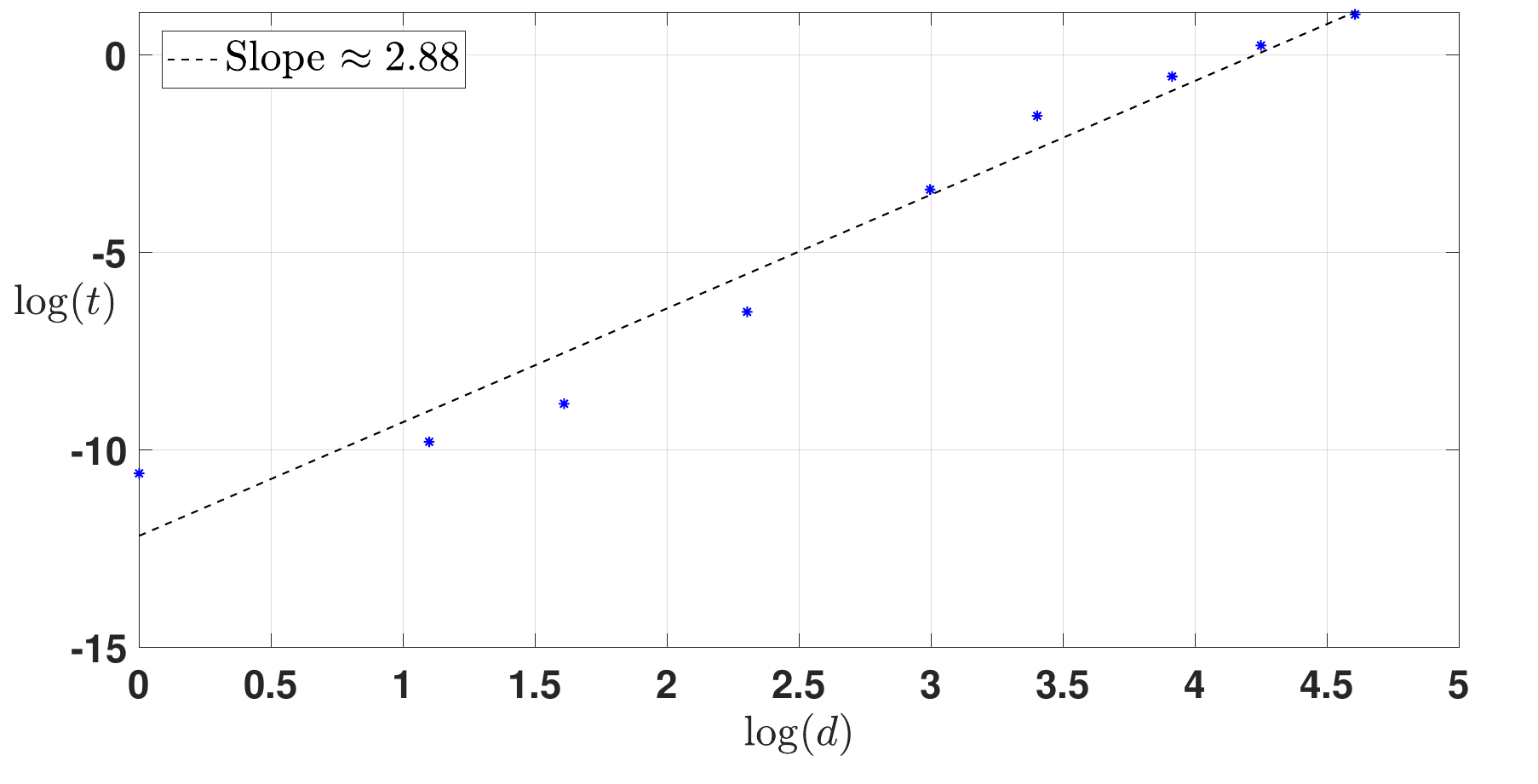}
	\caption{The relation between log CPU times and log state dimension $d$ with a polynomial slope of $\approx 2.88$ for the system using \eqref{eqn-4.19}.}
	\label{cpu_time}
\end{figure}

To examine the accuracy of the multivariate decomposition method with increasing dimensionality, the Mean Relative Error (MRE) across the same range of dimensions as in Fig. \ref{cpu_time} is presented in Table \ref{table-2}, where
\begin{equation*}
MRE:=\frac{\sum\limits_{k=0}^{N_t}\left|X_{t_k}-\hat X_{t_k}\right|}{\sum\limits_{k=0}^{N_t}|X_{t_k}|},
\end{equation*}
with $X_t$ and $\hat X_t$ denoting the true state and its estimate at time $t$, $|\cdot|$ representing the Euclidean norm, and $N_t$ being the total number of experimental time steps. Table \ref{table-2} shows that for a fixed number of particles $N_p=50$, the MREs remain bounded and at an acceptable level, exhibiting only slight growth as the dimension increases. This insensitivity to the state space dimension is a key indicator of the method’s robustness, enabling it to tackle the challenges inherent in high-dimensional scenarios.
\begin{table}[ht!]
    \centering
\begin{tabular}{|c|c|c|c|}
\hline
 Dimension [d] & $1$ & $3$ & $5$  \\
    \hline
    MRE & $0.4011$ & $0.4338$ & $0.4522$ \\
    \hline
     Dimension [d] & $10$& $20$& $30$ \\
     \hline
      MRE  & $0.4912$ & $0.5519$& $0.5568$ \\
     \hline
   Dimension [d]   & $50$ & $70$ & $100$ \\
    \hline
      MRE  & $0.5731$ & $0.5936$ & $0.6434$\\
      \hline
\end{tabular}
    \vspace{1em}\caption{The MREs versus the increasing dimension}
    \label{table-2}
\end{table}
\begin{figure}[ht!]
	\centering
        \includegraphics[width=0.5\textwidth, trim = 5 180 10 200,clip]{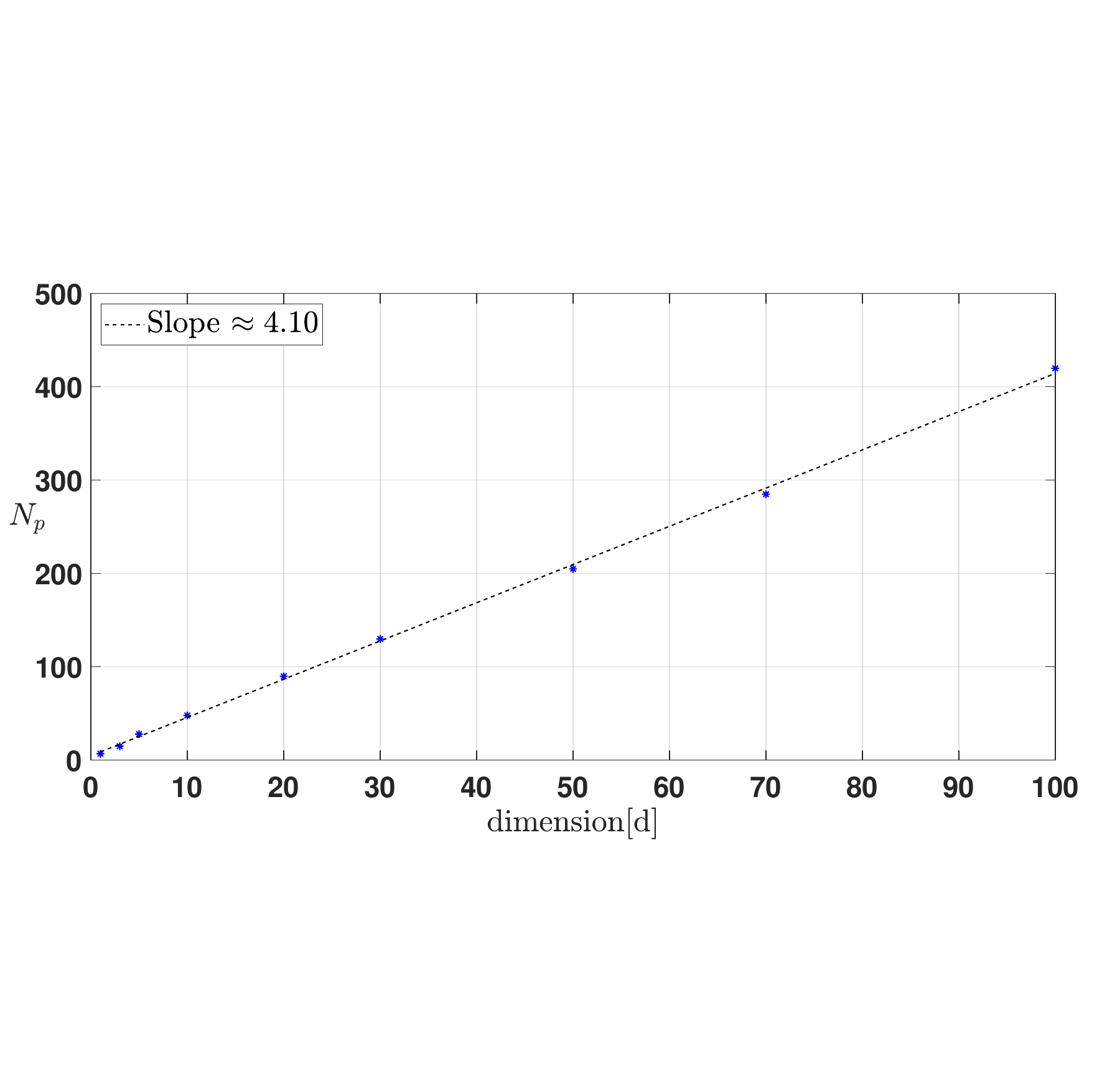}
	\caption{The relation between the number of particles and the dimension under the control of MRE $\leq 0.4$.}
	\label{dim_Np}
\end{figure}

To characterize the relationship between the number of particles and state dimension, we fix the MRE at $\leq0.4$ and record the required number of particles $N_p$ across the same dimension range as before, as shown in Fig. \ref{dim_Np}. The results indicate that the required $N_p$ exhibits approximately linear growth with increasing dimensionality. This behavior may stem from the decoupled and sparse structure of this specific example, which warrants further investigation. Nevertheless, this scaling trend suggests that the proposed multivariate decomposition method holds significant potential for mitigating the curse of dimensionality to a considerable extent.
\end{example}

\subsection{Numerics}\label{Numerics}

In this subsection, we shall numerically solve two benchmark multivariate NLF problems to illustrate the tracking capability and efficiency of our multivariate decomposition method for FPF. 

\begin{example}\label{ex-4.3}(The ship's tracking problem \cite{ALC:07})

This example describes that a ship moves with constant positive radial and angular velocities, though these are perturbed by white noise, when it is within a certain distance of the origin. Should it drift excessively far from the origin, a restoring force acts to push it back toward the origin. The signal model for the state process $X_t = \left(X_{t,1},X_{t,2}\right)^T\in\R^2$ is described by
\begin{equation}
\left\{
\begin{aligned}\label{eqn-4.2}
dX_{t,1}&=-X_{t,2}dt+f_1(X_{t,1},X_{t,2})dt+QdB_{t,1}\\
dX_{t,2}&=X_{t,1}dt+f_2(X_{t,1},X_{t,2})dt+QdB_{t,2}
\end{aligned}
\right.,
\end{equation}
where $B_{t}=\left(B_{t,1},B_{t,2}\right)^T$ is $2$-dimensional independent white noise process, and
\begin{equation}
f_i(\bx)=\gamma\frac{x_{j}}{|\bx|^2}-\theta\frac{x_{j}}{|\bx|}1_{(\rho,\infty)}(|\bx|), 
\end{equation}
$j=1,2$, where $|\cdot|$ is the Euclidean distance, $1_{(\rho,\infty)}(|\bx|)$ denotes the indicator function on the interval $(\rho,\infty)\subset\R$, and $\gamma$, $\theta$, $\rho$ are real-valued parameters. The observation is made every $\Delta t=0.05$ via the angular measurement:
\begin{equation}\label{eqn-4.1}
    Z_{t_k}=\arctan\left(\frac{X_{t_k,2}}{X_{t_k,1}}\right)+RW_{t_k},
\end{equation}
$k=1,\cdots,N_t$, where $N_t=\frac T{\Delta t}$, $W_{t_k}$ are i.i.d. $\N (0,1)$, independent of $(X_0,B_t)$. In the following simulations, we set $R=0.32$, $\gamma = 2$, $\theta = 50$, $\rho = 9$ and the initial condition $\bX_0=(0.5,-0.5)^T$, the same as those in \cite{ALC:07}. 
\end{example}

In order to apply the decomposition method to \eqref{eqn-4.2}-\eqref{eqn-4.1}, the observation function $h$ should be a polynomial in the state. It is natural to convert the system \eqref{eqn-4.2}-\eqref{eqn-4.1} from Cartesian coordinate into polar one. Let $X_{t,1}=\rho_t\cos\theta_t$, $X_{t,2}=\rho_t\sin\theta_t$, with $\rho_t\in\R_+$ and $\theta_t\in\R$, then
\begin{enumerate}
	\item The state process $\bx_t=(x_{t,1},x_{t,2}):=(\theta_t,\rho_t)$ is 
	\qquad \begin{equation}
	\left\{
	\begin{aligned}
	d\theta_t=&dt-Q\frac{\sin\theta_t}{\rho_t}dB_{t,1}+Q\frac{\cos\theta_t}{\rho_t}dB_{t,2} \\
	d\rho_t=&\frac{\gamma+\frac{Q^2}{2}}{\rho_t}dt-\theta1_{(\rho,\infty)}(\rho_t)dt+Q\cos\theta_tdB_{t,1}\\
    &+Q\sin\theta_tdB_{t,2}
	\end{aligned}
	\right.,
	\end{equation} 
    \item The observation process is
	\begin{equation}
	Z_{t_k}=\theta_{t_k}+RW_{t_k},
	\end{equation} 
    where $k=1,\cdots,N_t$.
\end{enumerate}

In this polar coordinate, the observation function $h(\bx)=x_1=\theta\in\R$, i.e. $m=1$. Let us assume that $\Sigma_i=\textup{diag}(\varepsilon_1,\varepsilon_2)$. Equation \eqref{eqn-7} becomes
\begin{align}\label{}\notag
	\left(\begin{array}{ccc}
	B_1&C_2&\\
	A_1&B_2&C_3\\
\end{array}\right)
       \ \left(\begin{array}{c}
	\tilde\Phi_{j,1}^i\\
	\tilde\Phi_{j,2}^i\\
        \tilde\Phi_{j,3}^i\\
\end{array}\right)
=
\left(\begin{array}{c}
	{\bf a}_{j,0}-C_j^i\\
	{\bf a}_{j,1}\\
\end{array}\right)&,
\end{align}
where
\begin{align*}
    \tilde\Phi_{j,1}^i=&\left(\tilde{\varphi}_{\be_1}^i,\tilde{\varphi}_{\be_2}^i\right)^T,\
    \tilde\Phi_{j,2}^i=\left(\tilde{\varphi}_{2\be_1}^i,\tilde{\varphi}_{\be_1+\be_2}^i,   \tilde{\varphi}_{2\be_2}^i\right)^T,\\
    \tilde\Phi_{j,3}^i=&\left(\tilde{\varphi}_{3\be_1}^i,\tilde{\varphi}_{2\be_1+\be_2}^i,
    \tilde{\varphi}_{\be_1+2\be_2}^i,
    \tilde{\varphi}_{3\be_2}^i\right)^T,\\
    B_1=&-\left(\frac{2X_{t,1}^i}{\varepsilon_1},\frac{2X_{t,2}^i}{\varepsilon_2}\right),\quad
A_1=\textup{diag}\left(\frac{1}{\varepsilon_1},\frac{1}{\varepsilon_2}\right),
\end{align*}
and
\begin{align*}
    {\bf a}_{j,0}=0,\quad{\bf a}_{j,1}=\left(\frac{1}{2},0\right)^T.
\end{align*}
It is clear to see that $A_1$ is invertible, which matches the conclusion in Proposition \ref{prop-3.6} or \ref{prop-3.7}. By backward recursion, we shall let $\tilde\Phi_{j,2}^i=\tilde\Phi_{j,3}^i=0$, so $B_2$, $C_2$ and $C_3$ have no impact on the coefficient $\tilde\Phi_{j,1}^i$,  then according to Proposition \ref{P2}, it yields that 
\begin{align*}
    \tilde\Phi_{j,1}^i=\left(\frac{\varepsilon_1}{2},    0\right)^T,\quad
C^i=X_{t,1}^i.
\end{align*}
Consequently,
\begin{equation*}
	\nabla\varphi^i(\bx)=(\varepsilon_{1},0)^T,
	\end{equation*} 
by \eqref{eqn-3.4}. By Proposition \ref{P1} and Theorem \ref{thm-3.5}, we have 
\begin{align*}
	\nabla\varphi(\bx)=&\left[\begin{array}{c}\varepsilon_{1}\\0\end{array}\right]
    +\frac{1}{\sum_{i=1}^{N_p}\mathcal{N}(\bx;X_t^i,\Sigma_i)}\sum_{i=1}^{N_p}
 (\bx-X_t^i)\\
    &\phantom{aa}\cdot\left(\frac1{N_p}\sum_{i=1}^{N_p}X_{t,1}^i-X_{t,1}^i\right)\frac{1}{2\pi|\Sigma_i|^{\frac{1}{2}}}\gamma\left(1,\frac{r^2}{2}\right)r^{-2}.
\end{align*}

\begin{figure}[htbp]
	\centering
        \includegraphics[width=0.5\textwidth, trim = 0 200 10 210,clip]{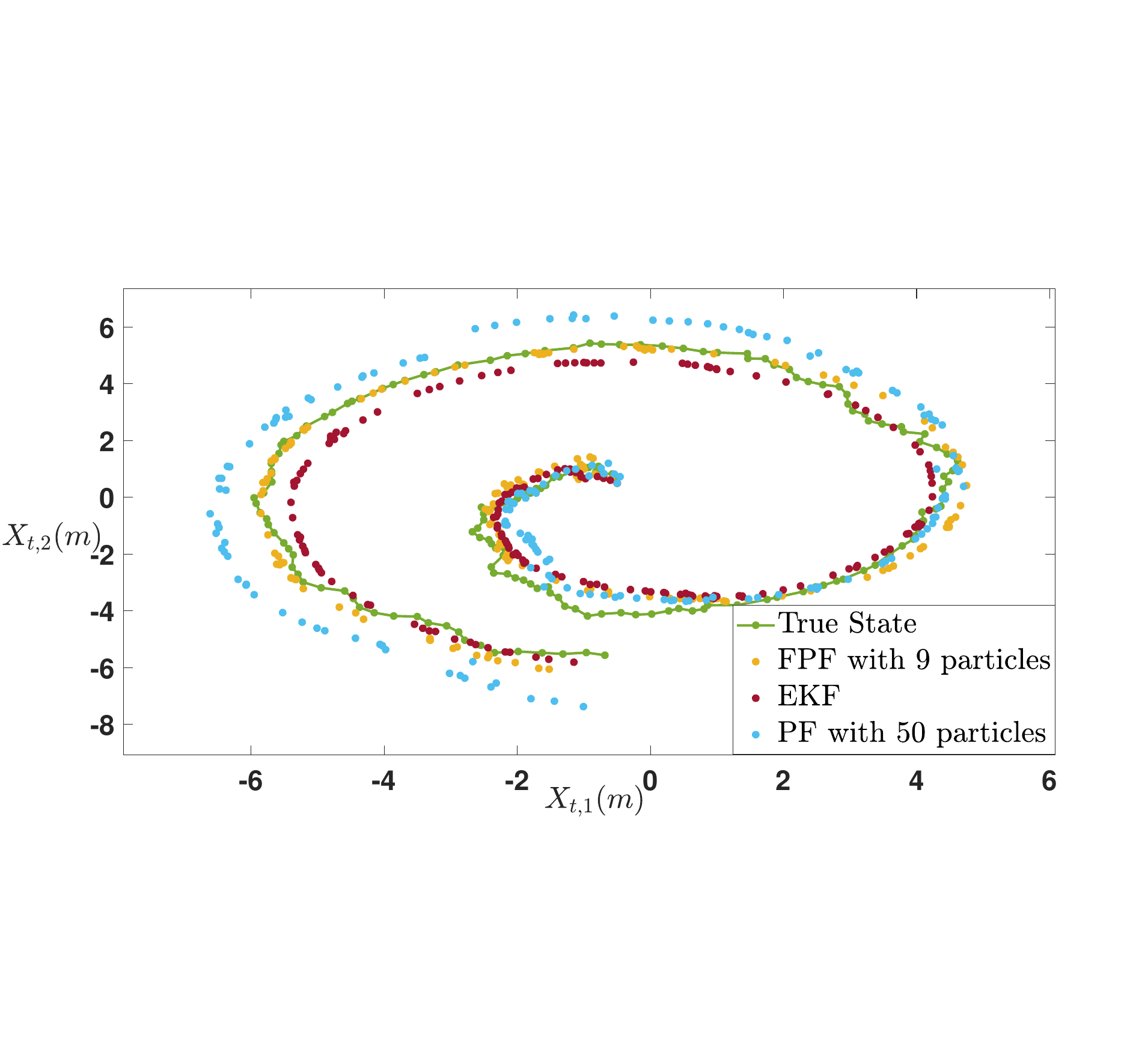}
	\caption{The tracking performances of the ship's trajectory with different methods. Each point represents the true/approximate location in the Cartesian coordinates at each time instant.}
	\label{d2_track}
\end{figure}

We next present the results of the numerical experiments. We set the total experimental time $T=8.25$, with time step $\Delta t = 0.05$. The covariance matrix $\Sigma_i=\textup{diag}(\varepsilon_1,\varepsilon_2)$ with $\varepsilon_1=\varepsilon_2=0.1$. Fig. \ref{d2_track} shows the tracking performance of FPF with $9$ particles, EKF and PF with $50$ particles.

Next, we do a Monte Carlo simulation with $M=100$ trajectories over the time interval $[0,T]$. The averaged root mean square error (ARMSE) over all trajectories and all time instants is defined as
\begin{equation}\label{eqn-4.23}
    \textup{ARMSE}^2=\frac{1}{M}\frac{1}{N_t}\sum_{j=1}^M\sum_{k=1}^{N_t}\left|X_{t_k,j}-\hat{X}_{t_k,j}\right|^2,
\end{equation}
where $X_{t_k,j}$ and $\hat{X}_{t_k,j}$ are the $j$th component of the true state and the estimated state at time instant $t_k$, respectively. The ARMSE and the CPU times of different methods are shown in Table \ref{table-1}. The FPF with $9$ particles has a smaller ARMSE than the PF with $50$ particles within $4$ times shorter CPU times. The FPF with less particle numbers beats the PF in both accuracy and efficiency. As for the EKF, it has a moderate ARMSE but the shortest CPU time, suggesting a trade-off between accuracy and efficiency.

\begin{table}[ht!]
    \centering
\begin{tabular}{|c|c|c|}
\hline
        & ARMSE & CPU times (s) \\
\hline
         FPF with $9$ particles & ${\bf 1.27506}$ & $0.0402099$\\
    \hline
    EKF & $1.38836$ & ${\bf 0.0062072}$\\
    \hline
    PF with $50$ particles & $1.68392$ & $0.1597206$\\
    \hline
\end{tabular}
\vspace{1em}    \caption{Comparison of FPF with $9$ particles, EKF and PF with $50$ particles in terms of ARMSE and CPU Times.}
\label{table-1}
    \end{table}

\begin{example}\label{ex-4.4}(The Lorenz oscillator system)

The Lorenz system, first developed based on the study of Rayleigh-B\'{e}nard convection, is a simple model that uses three interconnected nonlinear differential equations to simulate convection in the atmosphere. In a seminal work published in 1963, E. Lorenz demonstrated the intrinsic capacity of the system to manifest chaotic behavior. It has become a benchmark example in the NLF's literature to demonstrate the accuracy and efficiency of any novel algorithm. The state's system is described as follows \cite{NAP:11}:
\begin{equation}
\left\{
\begin{aligned}\label{NE4}
dX_{t,1}&=-\sigma(X_{t,1}-X_{t,2})dt+QdB_{t,1}\\
dX_{t,2}&=(-X_{t,1}X_{t,3}+\rho X_{t,1}-X_{t,2})dt+QdB_{t,2} \\
dX_{t,3}&=(X_{t,1}X_{t,2}-\beta X_{t,3})dt+QdB_{t,3}
\end{aligned}
\right.,
\end{equation}
where $Q=0.18$, $\sigma=10$, $\beta=\frac{8}{3}$ and $\rho= 25$. 
The observation process is
\begin{equation}\label{eqn-4.6}
    dZ_{t}=X_{t,1}dt+RdW_{t},
    \end{equation}
where $R=0.2$.
\end{example}

Similar as the ship's trajectory tracking problem in Example \ref{ex-4.3}, we suppose the variance matrix is $\Sigma_i=\textup{diag}(\varepsilon_1,\varepsilon_2,\varepsilon_3)$. Equation \eqref{eqn-7} for the Lorenz system \eqref{NE4} becomes
\begin{align}\label{}\notag
	\left(\begin{array}{ccc}
	B_1&C_2&\\
	A_1&B_2&C_3\\
\end{array}\right)
       \ \left(\begin{array}{c}
	\tilde\Phi_{j,1}^i\\
	\tilde\Phi_{j,2}^i\\
        \tilde\Phi_{j,3}^i\\
\end{array}\right)
=
\left(\begin{array}{c}
	{\bf a}_{j,0}-C_j^i\\
	{\bf a}_{j,1}\\
\end{array}\right),
\end{align}
where
\begin{align*}
    \tilde\Phi_{j,1}^i=&\left(\tilde{\varphi}_{\be_1}^i,
    \tilde{\varphi}_{\be_2}^i,
    \tilde{\varphi}_{\be_3}^i\right)^T,\\
    \tilde\Phi_{j,2}^i=&\left(\tilde{\varphi}_{2\be_1}^i,
    \tilde{\varphi}_{\be_1+\be_2}^i,
    \tilde{\varphi}_{\be_1+\be_3}^i,
    \tilde{\varphi}_{2\be_2}^i,
    \tilde{\varphi}_{\be_2+\be_3}^i,
    \tilde{\varphi}_{2\be_3}^i\right)^T,\\
    \tilde\Phi_{j,3}^i=&\left(\tilde{\varphi}_{3\be_1}^i,
    \tilde{\varphi}_{2\be_1+\be_2}^i,
    \tilde{\varphi}_{2\be_1+\be_3}^i,
    \tilde{\varphi}_{\be_1+2\be_2}^i,
    \tilde{\varphi}_{\be_1+2\be_3}^i,
    \tilde{\varphi}_{\be_1+\be_2+\be_3}^i,\right.\\
    &\phantom{aa}\left.\tilde{\varphi}_{3\be_2}^i,
    \tilde{\varphi}_{2\be_2+\be_3}^i,
    \tilde{\varphi}_{\be_2+2\be_3}^i,
    \tilde{\varphi}_{3\be_3}^i\right)^T,\\
    B_1=&-\left(\frac{2X_{t,1}^i}{\varepsilon_1},\frac{2X_{t,2}^i}{\varepsilon_2},\frac{2X_{t,3}^i}{\varepsilon_3}\right),\ 
    A_1=\textup{diag}\left(\frac{1}{\varepsilon_1} ,\frac{1}{\varepsilon_2},\frac{1}{\varepsilon_3}\right),
\end{align*}
and
\begin{align*}
    {\bf a}_{j,0}=0,\quad{\bf a}_{j,1}=\left(\frac{1}{2},0,0\right)^T.
\end{align*}
It is clear to see that $A_1$ is invertible, which matches the conclusion in Proposition \ref{prop-3.6} or \ref{prop-3.7}. Similar as in Example \ref{ex-4.3}, by backward recursion, letting $\tilde\Phi_{j,2}^i=\tilde\Phi_{j,3}^i=0$,  then according to Proposition \ref{P2}, it yields that 
\begin{align}\label{eqn-NS-3.3}
    \tilde\Phi_{j,1}^i=\left(\frac{\varepsilon_1}{2},0,0\right)^T,\quad
C^i=X_{t,1}^i.
\end{align}
thus,
\begin{equation*}
    \nabla\varphi^i(\bx)=\left(\varepsilon_{1},0,0\right)^T,
\end{equation*}
by \eqref{eqn-3.4}. By Proposition \ref{P1} and Theorem \ref{thm-3.5}, we have 
\begin{align*}
	&K_{3\times 1}(\bx)=\nabla\varphi(\bx)\\
    =&\left[\begin{array}{c}\varepsilon_{1}\\0\\0\end{array}\right]+\frac{1}{\sum_{i=1}^{N_p}\mathcal{N}(\bx;X_t^i,\Sigma_i)}\sum_{i=1}^{N_p}
 (\bx-X_t^i)\\
&\phantom{aaaaaaa}\cdot\left(\frac1{N_p}\sum_{i=1}^{N_p}X_{t,1}^i-X_{t,1}
^i\right)\frac{1}{2\pi^{\frac32}|\Sigma_i|^{\frac{1}{2}}}\gamma\left(\frac32,\frac{r^2}{2}\right)r^{-3}.
\end{align*}

\begin{figure}[ht!]
\setlength{\abovecaptionskip}{4pt}
\centering
\begin{minipage}[t]{0.5\textwidth}
    \centering 
    \includegraphics[width=0.9\linewidth, trim = 10 190 10 220,clip]{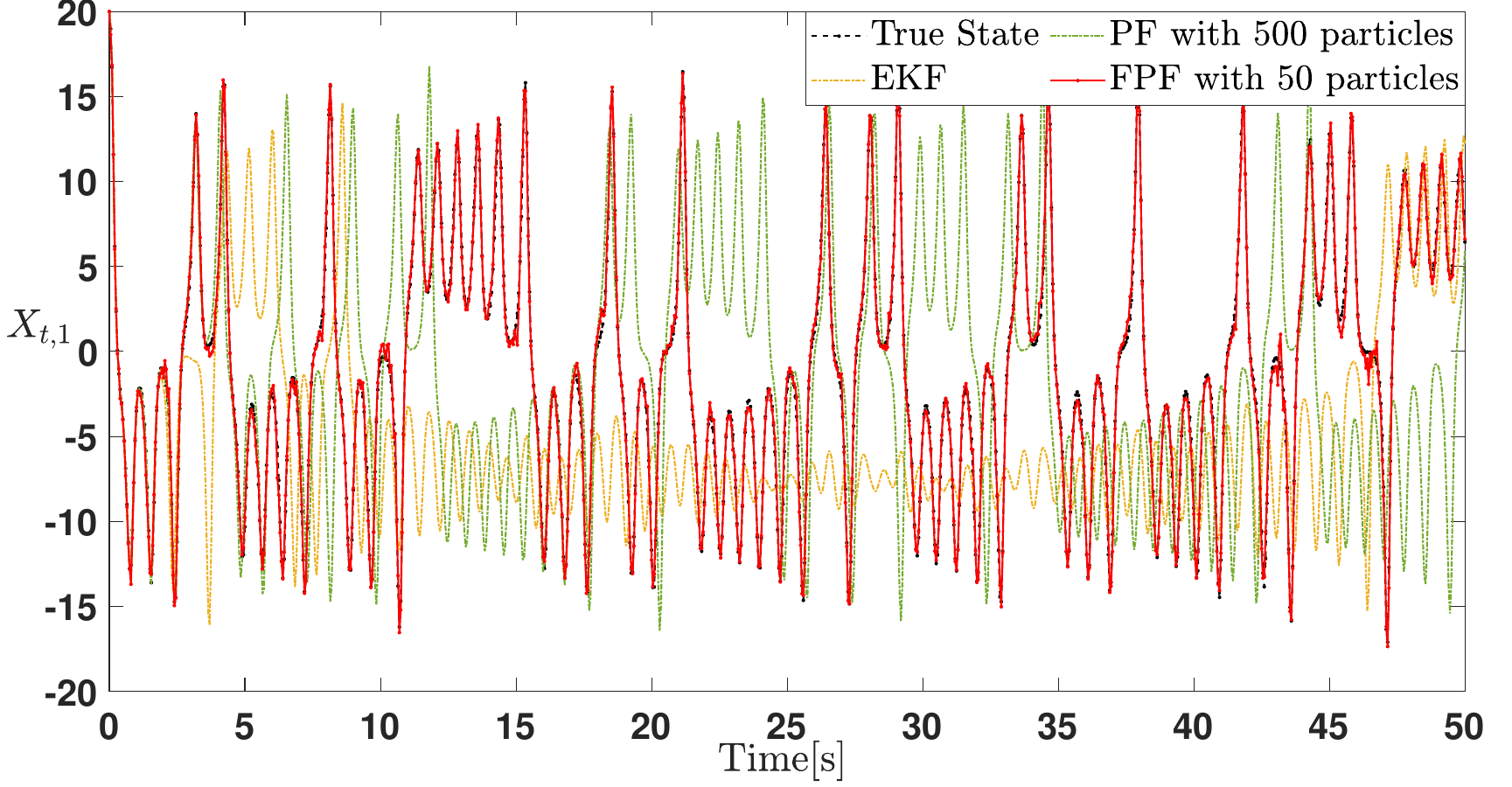}
    \small 
\end{minipage}
\vspace{4pt}

\begin{minipage}[t]{0.5\textwidth}
    \centering  
    \includegraphics[width=0.9\linewidth, trim = 10 190 10 220,clip]{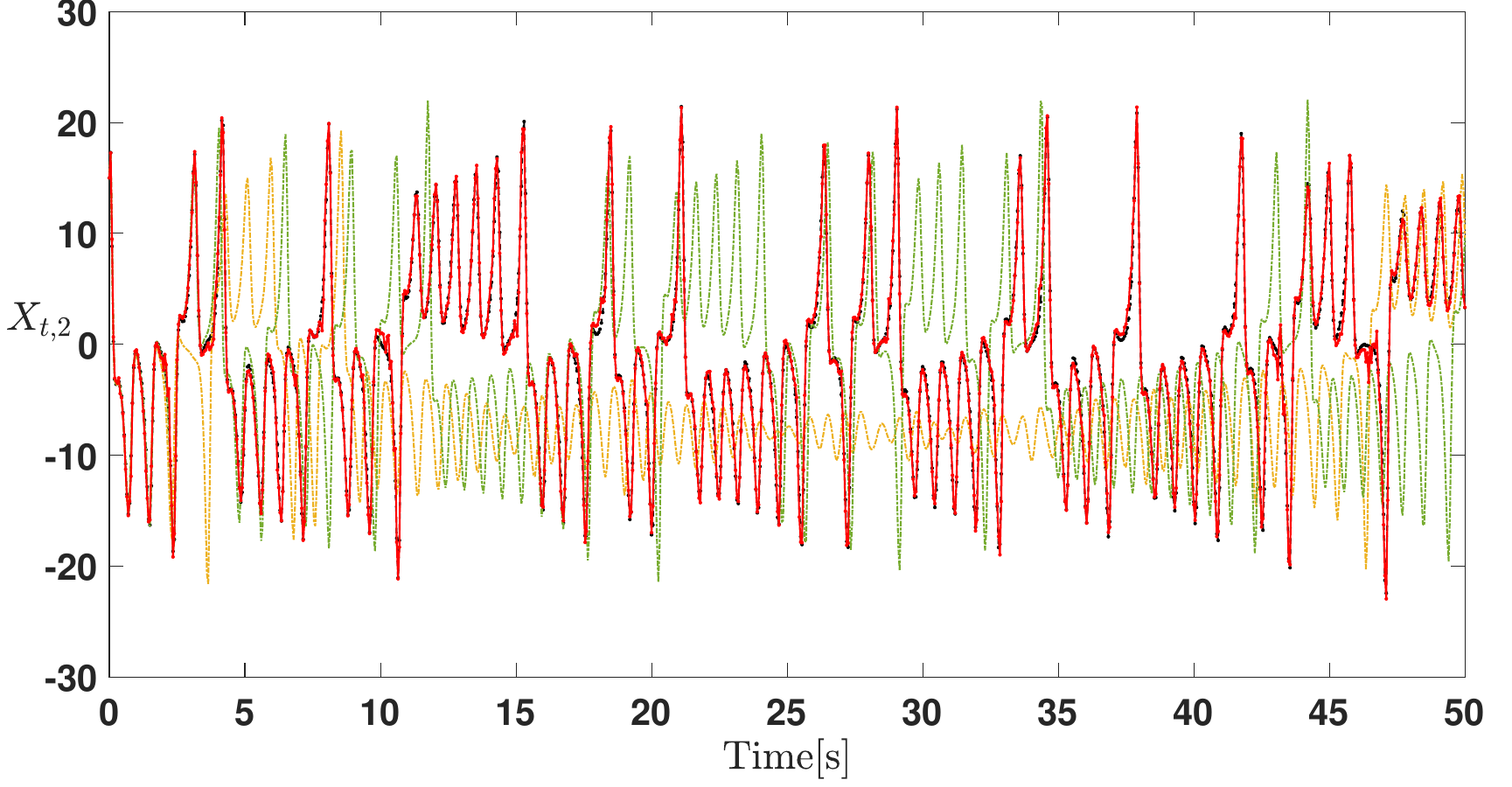}
    \small 
\end{minipage}
\vspace{4pt}

\begin{minipage}[t]{0.5\textwidth}
    \centering 
    \includegraphics[width=0.9\linewidth,trim = 10 190 10 220,clip]{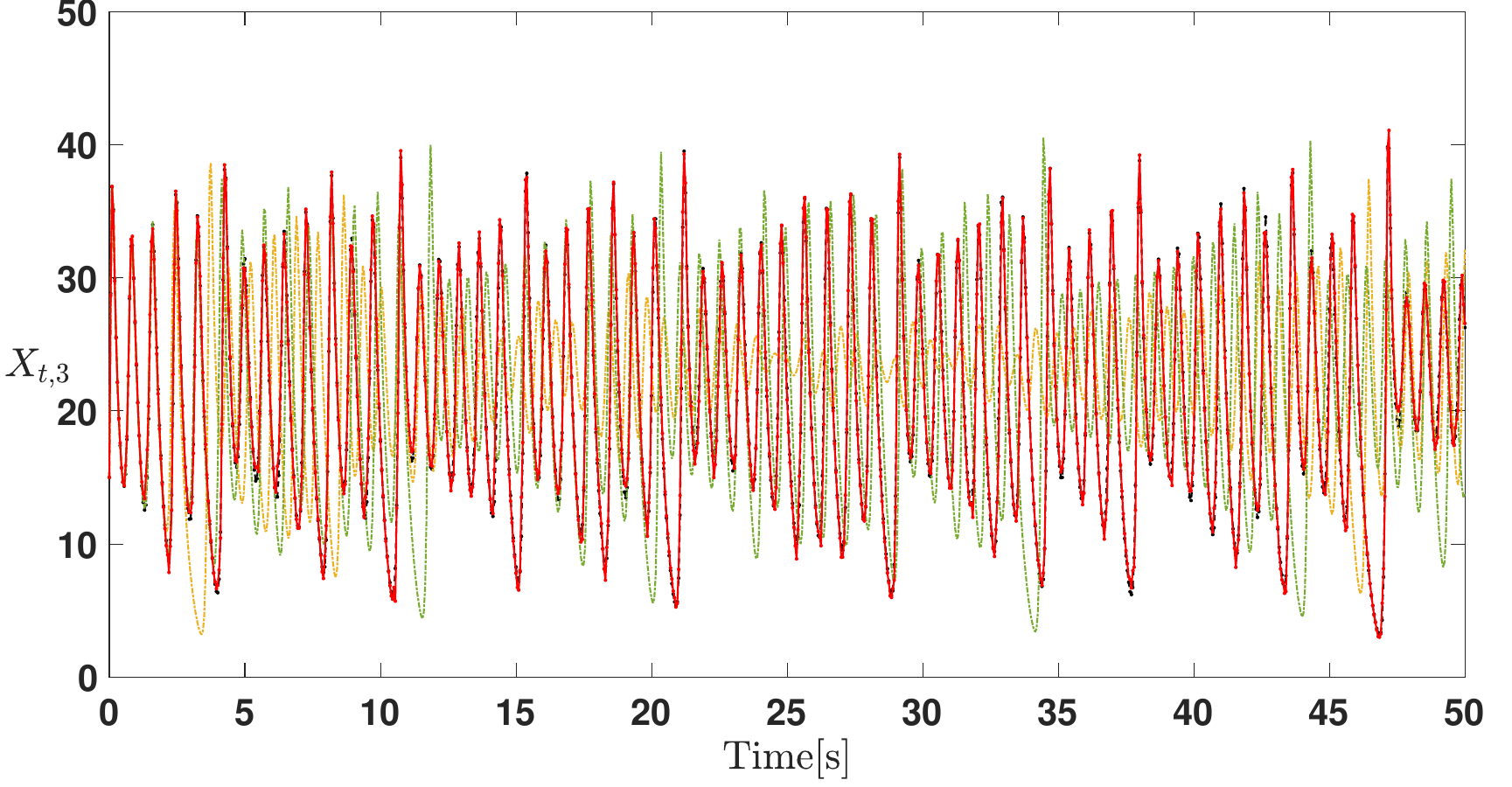}
    \small 
\end{minipage}

\vspace{4pt}
\small
\caption{The estimations of the true state (black) obtained by the FPF with the multivariate decomposition method (red), the EKF (yellow), and the PF (green), respectively.}
\label{lorenz_1}
\end{figure}

In the numerical simulation, we set the total experimental time $T=50$, with time step $\Delta t=0.001$. The true state is initialized as $X_0=(X_{0,1},X_{0,2},X_{0,3})^T=(20,15,15)^T$. The covariance matrix $\Sigma_i$ with $\varepsilon_1=\varepsilon_2=\varepsilon_3=0.01$. We use three different methods: the EKF, the PF with $500$ particles and the FPF with $50$ particles, to estimate the Lorenz system \eqref{NE4} and compare their performances in one realization in Fig. \ref{lorenz_1}. The CPU time of the FPF with the gain function obtained by our decomposition method is only $4.069$s, which is $9.38\times10^{-5}$s for each time step, while that of the PF is $26.700$s, that is $5.34\times10^{-4}$s for each time step.

\begin{table}[ht!]
    \centering
\begin{tabular}{|c|c|c|c|}
\hline
ARMSE & EKF  & PF with $500$ particles & FPF with $50$ particles  \\
    \hline
    $1$-$5$s  & $11.2471$ & $0.9956$ & ${\bf 0.7947}$ \\
    \hline
    $1$-$15$s  & $15.9726$ & $11.2724$ & ${\bf1.1211}$  \\
    \hline
    $1$-$30$s  & $16.4083$& $16.7982$& ${\bf1.1004}$ \\
    \hline
    $1$-$50$s  & $17.3346$& $18.1869$& ${\bf1.0879}$ \\
    \hline
\end{tabular}
    \vspace{1em}\caption{The ARMSE comparisons of the EKF, the PF with $500$ particles and the FPF with $50$ particles in different time periods.}
    \label{table-3}
\end{table}
To demonstrate the accuracy of the FPF, we calculate the ARMSE under the settings $M=1$, $T=50$, and $\Delta t=0.001$ across different time periods ($1$-$5$s, $1$-$15$s, $1$-$30$s, and $1$-$50$s), with results reported in Table \ref{table-3}. It can be observed that the FPF achieves the smallest ARMSE for all three estimated states across all time periods. The PF deteriorates rapidly from $8$s, which aligns with the tracking results of $X_{t,1}$ and $X_{t,2}$ in Fig. \ref{lorenz_1}. Regardless of the method employed, the error for $X_{t,1}$ is smaller than that for the other two states-this is likely because this state is directly observed, as per \eqref{eqn-4.6}. As the time period lengthens, the performance of the PF and EKF degrades, whereas that of the FPF remains around 1 throughout.

To further demonstrate the accuracy and efficiency of the FPF with different gain function approximations, we compare the multivariate decomposition method with the constant-gain approximation \cite{YLMM:13} and the kernel-based approach \cite{TM:16}, both of which are frequently used approximations in the literature of FPF. We set the total experimental time $T=10$, with time step $\Delta t=0.001$. The initial state $X_0$ and the covariance matrix $\Sigma_i$ are as before. In the kernel-based method, the bandwidth of kernel is set to be $\varepsilon = 0.1$ and the number of iterations to approximate the semi-group operator is $T_{iter} = 100$. 

\begin{figure}[ht!]
      \centering
      \includegraphics[width=0.5\textwidth, trim = 0 190 10 220,clip]{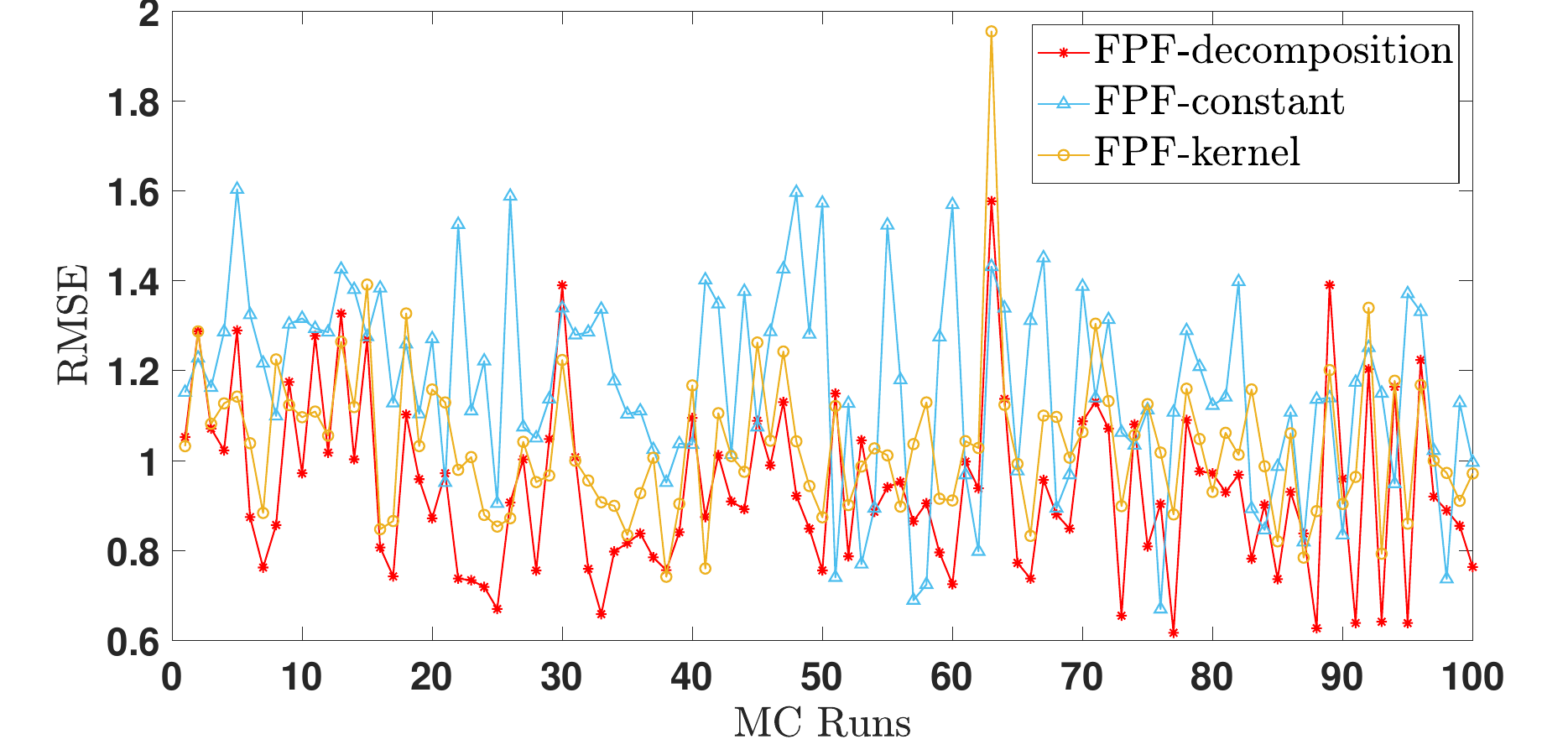}
    \caption{The RMSEs of every $100$ MC simulation of the FPF with three different gain function approximations.}\label{fig-6}
\end{figure}

The RMSEs of FPF with $50$ particles using different gain function approximations in $M=100$ MC simulations are plotted in Fig. \ref{fig-6}, where RMSE is the ARMSE defined in \eqref{eqn-4.23} with $M=1$. It shows that in nearly all simulations the multivariate decomposition method yields smaller RMSE than the other two, except a few. 

The ARMSE over $100$ MC simulations in \eqref{eqn-4.23} and the CPU times of those results in Fig. \ref{fig-6} are list in Table \ref{table-4}. In terms of ARMSE, the multivariate decomposition method outperforms both the constant-gain approximation and the kernel-based approximation, indicating higher estimation accuracy. For CPU times, the constant-gain approximation is the fastest, while the multivariate decomposition method is much more efficient than the kernel-based approximation, making the former a favorable balance of accuracy and computational efficiency.
\begin{table}[ht!]
    \centering
\begin{tabular}{|c|c|c|}
\hline
        & ARMSE & CPU times (s) \\
\hline
         multivariate decomposition method & ${\bf 0.9374}$ & $1.53$\\
    \hline
    constant-gain approximation & $1.1660$ & ${\bf 0.51}$\\
    \hline
    kernel-based approximation & $1.0348$ & $12.35$\\
    \hline
\end{tabular}
\vspace{1em}    \caption{Comparison of the FPF with $50$ particles using different gain function approximations in terms of ARMSE and CPU Times.}
\label{table-4}
    \end{table}


Finally, we briefly discuss the covariance matrix $\Sigma_i$ in the multivariate decomposition method. For simplicity and to ensure invertibility, we assume $\Sigma_i$ is a diagonal matrix (i.e., $\Sigma_i = \textup{diag}(\varepsilon 1_{d \times 1})$); however, the value of $\varepsilon$ significantly impacts the method’s performance. Here, we conduct experiments on the Lorenz system \eqref{NE4} under the same settings as before: $T=50$, $\Delta t=0.001$, and $N_p=50$. The RMSEs of each component $X_{t,j}$, $j=1,2,3$, for different $\varepsilon$ values are presented in Table \ref{table-5}. The results show that smaller $\varepsilon$ values tend to yield better estimation performance, as reflected by the monotonic decrease in RMSE for all three $X_{t,j}$s. Nevertheless, the rate of performance improvement slows as $\varepsilon$ decreases-this implies the existence of a diminishing returns point, beyond which further reducing $\varepsilon$ provides limited practical benefits.
\begin{table}[ht!]
    \centering
\begin{tabular}{|c|c|c|c|}
\hline
    \multirow{2}{*}{} & \multicolumn{3}{c|}{RMSE} \\ \cline{2-4}
         $\varepsilon$ & $X_1$ & $X_2$ & $X_3$ \\
    \hline
     $3$ & $0.8305$ & $1.6654$ & $1.5666$ \\
     \hline
     $0.3$ & $0.6498$ & $1.0882$ & $0.9430$ \\
     \hline
     $0.03$ & $0.6072$ & $0.9975$ & $0.8926$ \\
    \hline
\end{tabular}
\vspace{1em}
    \caption{Comparison of different $\varepsilon$ and RMSEs.}
    \label{table-5}
\end{table}

The selection of $\varepsilon$ presents a challenge in identifying a value that performs robustly across all problem types. While this issue has been explored in several studies, no practical guidelines have yet been established.
\begin{itemize}
        \item Silverman \cite{KS:87} proposed a parameter selection rule for $\varepsilon$, where $\varepsilon \propto N_p^{-\frac{1}{d+4}}$, where $d$ is the state's dimension.
        \item In the previous work of A. Taghvaei et. al\cite{TMM:20}, they highlights the presence of an optimal $\varepsilon$ for error minimization. Theoretically, this value balances the bias-variance trade-off, with the constant-gain case corresponding to $\varepsilon\to\infty$. An analytical formula for the optimum remains elusive. Instead, the heuristic $\varepsilon=\frac{4(\mathrm{med})^2}{\log(N)}$ is frequently used, where \textup{med} is the median pairwise distance \cite{ADCLS:17}. This choice maintains the kernel matrix sufficiently distinct from the identity matrix, thereby preventing degeneracy.
    \end{itemize}

\section{Conclusions}

\IEEEPARstart{T}{his} paper extends the decomposition method for the feedback particle filter (FPF) from one-dimensional to multivariate settings. A key challenge in this extension lies in the fact that the scalar systems’ analytical tractability vanishes in multivariate cases, as Poisson’s equation involves coupled terms and the divergence operator introduces non-trivial boundary constraints. Building on our prior work—where Poisson’s equation was split into two exactly solvable sub-equations for scalar systems with polynomial observations—we address these hurdles by innovatively adopting tensor product Hermite polynomials in the Galerkin method in the backward recursion for one sub-equation and the construction of the weighted radially symmetric solution to the other one. The main contributions include: resolving the sub-equations’ explicit solvability via this tailored spectral approach; rigorously proving the invertibility of coefficient matrices in some typical cases; and developing a unified framework that retains polynomial complexity. This yields an implementable multivariate FPF framework, formalized in Theorem \ref{thm-3.5}, which provides efficient, accurate approximate gain functions for polynomial observation systems.

To validate the gain function computed by the multivariate decomposition method, we compare it with the constant-gain approximation and kernel-based approach against the true gain function in a two-dimensional setting, with the decomposition method achieving the smallest $l^2$-error. We further apply the method to high-dimensional nonlinear systems, examining how CPU runtime scales with dimensionality and the number of particles needed to maintain accuracy. Comprehensive comparisons with established methods, the EKF, PF, kernel-based FPF, and constant-gain FPF, across two practical numerical examples demonstrate the decomposition method’s strong robustness and superior computational accuracy. Moreover, it significantly outperforms the PF and kernel-based approach in computational efficiency.

\bibliographystyle{plain} 
\bibliography{references}     

\begin{thebibliography}{10}

\bibitem{CASV:04}
C.~Andrieu, A.~Doucet, S.S. Singh, and V.B. Tadic.
\newblock Particle methods for change detection, system identification, and
  control.
\newblock {\em Proceedings of the IEEE}, 92(3):423--438, 2004.

\bibitem{AMGC:02}
M.~Arulampalam, S.~Maskell, N.~Gordon, and T.~Clapp.
\newblock A tutorial on particle filters for online nonlinear/non-{G}aussian
  {B}ayesian tracking.
\newblock {\em IEEE Transactions on Signal Processing}, 50(2):174--188, 2002.

\bibitem{BLK:01}
Y.~Bar-Shalom, X.~Li, and T.~Kirubarajan.
\newblock {\em Estimation with Applications to Tracking and Navigation}.
\newblock New York : Wiley, 2001.

\bibitem{B:18}
K.~Berntorp.
\newblock Comparison of gain function approximation methods in the feedback
  particle filter.
\newblock In {\em In 2018 21th International Conference on Information Fusion},
  pages 123--130, 2018.

\bibitem{BG:16}
K.~Berntorp and P.~Grover.
\newblock Data-driven gain computation in the feedback particle filter.
\newblock In {\em Proceedings of the American Control Conference}, pages
  2711--2716, 2016.

\bibitem{ALC:07}
A.~Budhiraja, L.~Chen, and C.~Lee.
\newblock A survey of numerical methods for nonlinear filtering problems.
\newblock {\em Physica D: Nonlinear Phenomena}, 230(1):27--36, 2007.

\bibitem{ADCLS:17}
A.~Chaudhuri, D.~Kakde, C.~Sadek, L.~Gonzalez, and S.~Kong.
\newblock The mean and median criteria for kernel bandwidth selection for
  support vector data description.
\newblock In {\em 2017 IEEE International Conference on Data Mining Workshops},
  pages 842--849, 2017.

\bibitem{E:94}
G.~Evensen.
\newblock Sequential data assimilation with a nonlinear quasi-geostrophic model
  using monte carlo methods to forecast error statistics.
\newblock {\em Journal of Geophysical Research}, 99(5):10143--10162, 1994.

\bibitem{EG:03}
G.~Evensen.
\newblock The ensemble \uppercase{K}alman filter: theoretical formulation and
  practical implementation.
\newblock {\em Ocean Dynamics}, 53(4):343--367, 2003.

\bibitem{GSS:93}
N.~Gordon, D.~Salmond, and A.~Smith.
\newblock Novel approach to nonlinear/non-{G}aussian {B}ayesian state
  estimation.
\newblock {\em IEE proceedings F. Communications, Radar and Signal Processing},
  140(2):107--113, 1993.

\bibitem{JA:70}
A.~Jazwinski.
\newblock Stochastic processes and filtering theory.
\newblock {\em Academic Press}, 1970.

\bibitem{K:60}
R.~Kalman.
\newblock A new approach to linear filtering and prediction problems.
\newblock {\em Transactions of the ASME. Series D. Journal of Basic
  Engineering}, 82(1):35--45, 1960.

\bibitem{KB:61}
R.~Kalman and R.~Bucy.
\newblock New results in linear filtering and prediction theory.
\newblock {\em Transactions of the ASME. Series D. Journal of Basic
  Engineering}, 83(1):95--108, 1961.

\bibitem{K:90}
H.~Kunita.
\newblock {\em Stochastic flows and stochastic differential equations},
  volume~24 of {\em Cambridge Studies in Advanced Mathematics}.
\newblock Cambridge University Press, Cambridge, 1990.

\bibitem{NAP:11}
K.~Nosrati, A.~Shokouhi, N.~Pariz, and A.~Azemi.
\newblock Chaos synchronization of fractional-order lorenz system with
  unscented {K}alman filter.
\newblock In {\em 2011 19th Iranian Conference on Electrical Engineering},
  2011.

\bibitem{AS:18}
A.~Radhakrishnan and S.~Meyn.
\newblock Feedback particle filter design using a differential-loss reproducing
  kernel \uppercase{H}ilbert space.
\newblock In {\em Proceedings of the American Control Conference}, pages
  329--336, 2018.

\bibitem{RM:19}
A.~Radhakrishnan and S.~Meyn.
\newblock Gain function tracking in the feedback particle filter.
\newblock In {\em Proceedings of the American Control Conference}, pages
  5352--5359, 2019.

\bibitem{BSN:04}
B.~Ristic and S.~Arulampalam.
\newblock {\em Beyond the Kalman Filter: Particle Filters for Tracking
  Applications}.
\newblock Boston, MA : Artech House, 2004.

\bibitem{SS:66}
S.~Schmidt.
\newblock Application of state-space methods to navigation problems.
\newblock {\em Advances in Control Systems}, 3(2):293--340, 1966.

\bibitem{KS:87}
B.~Silverman.
\newblock {\em Density Estimation for Statistics and Data Analysis}.
\newblock Boca Raton: CRC Press, 1986.

\bibitem{SKP:19}
S.~Surace, A.~Kutschireiter, and J.-P. Pfister.
\newblock How to avoid the curse of dimensionality: scalability of particle
  filters with and without importance weights.
\newblock {\em SIAM Review}, 61(1):79--91, 2019.

\bibitem{TM:16}
A.~Taghvaei and P.~Mehta.
\newblock Gain function approximation in the feedback particle filter.
\newblock In {\em Proceedings of 2016 IEEE Conference on Decision and Control},
  pages 5446--5452, 2016.

\bibitem{AP:23}
A.~Taghvaei and P.~Mehta.
\newblock A survey of feedback particle filter and related controlled
  interacting particle systems ({CIPS}).
\newblock {\em Annual Reviews in Control}, 55:356--378, 2023.

\bibitem{TMM:17}
A.~Taghvaei, P.~Mehta, and S.~Meyn.
\newblock Error estimates for the kernel gain function approximatioin in the
  feedback particle filter.
\newblock In {\em Proceedings of the American Control Conference}, pages
  4576--4582, 2017.

\bibitem{TMM:20}
A.~Taghvaei, P.~Mehta, and S.~Meyn.
\newblock Diffusion map-based algorithm for gain function approximation in the
  feedback particle filter.
\newblock {\em SIAM/ASA Journal on Uncertainty Quantification},
  8(3):1090--1117, 2020.

\bibitem{WML:25}
R.~Wang, H.~Miao, and X.~Luo.
\newblock A decomposition approach for the gain function in the feedback
  particle filter.
\newblock accepted by Proceedings of 2025 IEEE Conference on Decision and
  Control, arXiv:2503.23662, 2025.

\bibitem{Y:14}
T.~Yang.
\newblock {\em Feedback particle filter and its applications}.
\newblock PhD thesis, University of Illinois at Urbana-Champaign, 2014.

\bibitem{YLMM:13}
T.~Yang, R.~Laugesen, P.~Mehta, and S.~Meyn.
\newblock Multivariable feedback particle filter.
\newblock In {\em Proceedings of 2012 IEEE Conference on Decision and Control},
  pages 4063--4070, 2012.

\bibitem{YLMM:16}
T.~Yang, R.~Laugesen, P.~Mehta, and S.~Meyn.
\newblock Multivariable feedback particle filter.
\newblock {\em Automatica}, 71:10--23, 2016.

\bibitem{YMM:13}
T.~Yang, P.~Mehta, and S.~Meyn.
\newblock Feedback particle filter.
\newblock {\em IEEE Transactions on Automatic Control}, 58(10):2465--2480,
  2013.

\bibitem{YMMS:12}
H.~Yin, P.~Mehta, S.~Meyn, and U.~Shanbhag.
\newblock Synchronization of coupled oscillators is a game.
\newblock {\em IEEE Transactions on Automatic Control}, 57(4):920--935, 2012.

\end{thebibliography}




\end{document}